\documentclass[12pt,letterpaper]{article}
\usepackage[latin1]{inputenc}
\usepackage{amsfonts}
\usepackage{amsmath}
\usepackage{amstext}
\usepackage{amsbsy}
\usepackage{amssymb}
\usepackage{graphics}
\usepackage{graphicx}
\usepackage{yfonts}
\usepackage{pstricks-add}
\usepackage{pst-math}
\usepackage{url}

\usepackage{authblk}

\newtheorem{thm}{Theorem}[section]
\newcommand{\bthm}{\begin{thm}\noindent{\bf }~}
\newcommand{\ethm}{\end{thm}}
\newtheorem{propo}{Proposition}[section]
\newcommand{\bprop}{\begin{propo}\noindent{\bf }~}
\newcommand{\eprop}{\end{propo}}
\newtheorem{lema}{Lemma}[section]
\newcommand{\blem}{\begin{lema}\noindent{\bf }~}
\newcommand{\elem}{\end{lema}}
\newtheorem{defe}{Definition}[section]
\newcommand{\bdefe}{\begin{defe}\noindent{\bf }~}
\newcommand{\edefe}{\end{defe}}
\newtheorem{ejem}{Example}[section]
\newcommand{\bejem}{\begin{ejem}\noindent{\bf }~}
\newcommand{\eejem}{\end{ejem}}
\newtheorem{coro}{Corollary}[section]
\newcommand{\bcor}{\begin{coro}\noindent{\bf }~}
\newcommand{\ecor}{\end{coro}}
\newtheorem{rema}{Remark}[section]
\newcommand{\brem}{\begin{rema}\noindent{\bf }~\rm}
\newcommand{\erem}{\end{rema}}
\newcommand{\bdemo}{\noindent\textbf{Proof. }~\rm}
\newcommand{\edemo}{\hfill{$\sqcap \kern-18pt \sqcup$}}

\title{\large\bf Pre-twisted algebrizable differential equations}

\author[D. Kopta et al.]
       {\small {Elifalet L\'opez-Gonz\'alez}
       \\ {\footnotesize{DM de la UACJ en Cuauhtémoc, Universidad Autónoma de Ciudad Juárez,
       Carretera Cuauhtémoc-Anáhuac, Km 3.5 S/N, Ejido Cuauhtémoc, C.P. 31600, Cd. Cuauhtémoc, Chih, México
       \ttfamily\\ elgonzal@uacj.mx
       }}}

\begin{document}

\maketitle

\begin{small}\noindent \textbf{Abstract.} We introduce the \emph{$\varphi\mathbb{A}$-differentiability} for
functions $f:U\subset \mathbb R^{k}\to\mathbb R^{n}$ where
$\mathbb A$ is the linear space $\mathbb R^{n}$ endowed with an
algebra product which is unital, associative, commutative, $U$ is
an open set, and $\varphi:U\subset \mathbb R^{k}\to\mathbb A$ is a
differentiable function in the usual sense. We call to this
differentiability \emph{pre-twisted differentiability}. We also
introduce the corresponding generalized Cauchy-Riemann equations
($\varphi\mathbb{A}$-CREs), the Cauchy-integral theorem, and the
\emph{$\varphi\mathbb A$-differential equations}. The
four-dimensional vector fields associated with triangular
billiards are $\varphi\mathbb{A}$-differentiable. The
\emph{$\varphi\mathbb A$-differential equations} can be used for
constructing exact solutions of partial differential equations of
mathematical physics like the three-dimensional heat, Laplace's,
and wave equations.
\end{small}

\noindent \textbf{Keyword}: \texttt{Partial differential
equations}, \texttt{Generalized Cauchy-Riemann equations},\\
\texttt{Lorch differentiability}



\section*{Introduction}

The theory of analytic functions over algebras has been developed
since the end of the 19th century, see \cite{AMEE}, \cite{JSC},
\cite{Ket}, \cite{Lor}, \cite{War1}, and \cite{War1}. There have
been several attempts to construct solutions of classical partial
differential equations (PDEs for plural and PDE for singular) of
mathematical physics by means of component functions of
differentiable functions in the sense of Lorch, also see
\cite{Pla}, \cite{Wag}, and \cite{War2}. Ketchum's work stands
out, which perhaps has not been understood or has been
misinterpreted, because he works with algebras in a general way,
which makes it seem a bit complicated to construct solutions of
PDEs by the method he proposes. In \cite{EL2} and \cite{EL3}
algebras have been given for which Ketchum's results are
simplified, since in these the conditions (2), (3), and (4), given
in \cite{Ket} pp. 642 are satisfied, so that their fulfillment is
no longer required. Thus, only the conditions given in (47) of pp.
660 are required to be satisfied, and if we consider $w=w(x)$
affine, only the second condition of (47) will be required.

In this paper a concept of differentiability is introduced; the
``Pre-twisted Differentiability'', which is similar to the Lorch
derivative. This differentiability can be used to make more
explicit the method given by P. W. Ketchum in \cite{Ket} for
constructing solutions of PDEs of mathematical physics. Therefore,
this new notion of differentiability constitutes an important key
to the construction of solutions of classical mathematical physics
PDEs, see \cite{EL2} and \cite{EL3}. In particular, families of
pre-twisted differentiable functions are bi-harmonic functions.

In \cite{Ket} pp. 659, 660, analytic functions $f(w)$ where
$w=w(x)$ and $f$ is differentiable en the sense of Lorch, are
considered. An algebra $\mathbb A$ whose analytic functions $f(w)$
satisfy the Laplace's equation is called \emph{harmonic algebra}.
In the paper, pp. 547, it is interpreted that for an algebra to be
a harmonic algebra it is required that $e_1^2+e_2^2+ e_3^2 =0$,
but this condition is necessary for the case of $w=(x,y,z)$ (in
our notation $\varphi(x,y,z)=(x,y,z)$). P. W. Ketchum makes a
similar claim in the introduction to \cite{Ket1}. This condition
is used for solving PDEs, see \cite{Pog}. In \cite{EL3} it is
shown that for the algebra $\mathbb A$ defined by $\mathbb R^3$
with the product $e_3^2=e_2$, $e_3^2=e_1$, where the unit $e$ of
the algebra is $e=e_1$, is a harmonic algebra by taking
$$
w=\varphi(x,y,z)=\left(x+k_1,-\frac{1}{2}x+\frac{\sqrt{2}}{2}y+\frac{1}{2}z+k_2,-\frac{1}{2}x-\frac{\sqrt{2}}{2}y-\frac{1}{2}z+k_3\right).
$$
With respect to $\mathbb A$ one has that $e_1^2+e_2^2+
e_3^2=e_1+e_3+e_1=(2,0,1)$. So, $e_1^2+e_2^2+e_3^2\neq 0$.

The pre-twisted differentiability depends on a differentiable
function in the usual sense $\varphi$ and on an unital associative
commutative algebra $\mathbb A$, see Section \ref{spadr}. Thus, to
make this dependence explicit, we call this $\varphi\mathbb
A$-differentiability. Tentatively, one can say that a function $f$
is $\varphi\mathbb A$-differentiable if $f=L\circ\varphi$ and its
derivative $f'(x_0)=L'\circ\varphi(x_0)$, where L is
differentiable in the sense of Lorch with respect to $\mathbb A$
and $L'$ is its derivative in the sense of Lorch. If
$\varphi(x,y)=(x,y)$ and $\mathbb A$ is the algebra defined by the
complex filed $\mathbb C$, then $\varphi\mathbb
A$-differentiability is the usual differentiability in the complex
sense.

The pre-twisted differentiability has associated systems of linear
partial differential equations, which under certain conditions
characterize this differentiability. These PDE systems can be in
the form
\begin{equation}\label{s1}
    \begin{array}{ccc}
  a_{111}u_x+a_{121}v_x+a_{131}w_x+a_{112}u_y+a_{122}v_y+a_{132}w_y  & = & 0 \\
  a_{211}u_x+a_{221}v_x+a_{231}w_x+a_{212}u_y+a_{222}v_y+a_{232}w_y  & = & 0 \\
  a_{311}u_x+a_{321}v_x+a_{331}w_x+a_{312}u_y+a_{322}v_y+a_{332}w_y  & = & 0 \\
\end{array},
\end{equation}
where $a_{ijk}$ are functions of $(x,y,z)$, $u_x=\frac{\partial
u}{\partial x}$ and so on. Therefore, these systems can be called
``\emph{Generalized Cauchy-Riemann equations}'', see (\ref{s1})
for the case of differentiability in the sense of Lorch.

The algebrizability of vector fields and differential equations
has been studied in \cite{AMEE}, \cite{AMEENADE}, \cite{AFLY},
\cite{Dyg}, \cite{Dyg3}, and \cite{EL1}. The case of
algebrizability of systems of two first order partial differential
equations with two dependent and two independent variables is
being developed in \cite{CAL1} and \cite{CAL2}. In this paper we
define the $\varphi\mathbb A$-algebrizability of vector fields and
autonomous ordinary differential equations. The corresponding
Cauchy-Riemann equations ($\varphi\mathbb A$-CREs) are introduced
and a version of the Cauchy-integral Theorem is showed. We give
examples of generalized systems of Cauchy-Riemann equations for
the $\varphi\mathbb A$-differentiability. The following inverse
problem is considered: when a given homogeneous linear system of
two first order partial differential equations corresponds to the
CREs for some function $\varphi$ and a two dimensional algebra
$\mathbb{A}$. We show that two-dimensional quadratic complex
vector fields associated with triangular billiards are
$\varphi\mathbb A$-algebrizable. For adequate algebra $\mathbb A$
and a function $\varphi$ we show that the set of $\varphi\mathbb
A$-differential functions are solutions of homogeneous first order
linear partial differential equations with two dependent and two
independent variables. In \cite{CAL1} and \cite{CAL2} solutions
are constructed for first order linear systems of two partial
differential equations with constant coefficients, see Section
\ref{algsection}.

The organization of this paper is the following. In Section
\ref{spadr} we recall the definition of an algebra which we denote
by $\mathbb A$, we introduce the $\varphi\mathbb
A$-differentiability and give some results related to this like
the $\varphi\mathbb A$-CREs. In Section \ref{s2} the problem of
when a given first order linear system of two partial differential
equations results the $\varphi\mathbb A$-CREs for some
differentiable function $\varphi$ and an algebra $\mathbb A$, is
considered. In Section \ref{s3} the $\varphi\mathbb A$-line
integral is defined, a corresponding Cauchy-integral theorem is
given, a generalization of the fundamental theorem of calculus is
a consequence obtained, and examples are given. Moreover, the
$\varphi\mathbb A$-differential equations are introduced, an
existence and uniqueness of solutions theorem is given, it is
showed that systems associated to triangular billiards are
$\varphi\mathbb A$-differentiable, and $\varphi\mathbb
M$-algebrizable differential equations are introduced for
commutative matrix algebras $\mathbb M$. In Section
\ref{sol:singles:pdes} we give examples of solutions of PDEs and
systems PDEs constructed by using $\varphi\mathbb
A$-differentiability and $\varphi\mathbb A$-differential
equations.

\section{$\varphi\mathbb A$-differentiability}\label{spadr}

\subsection{Algebras}

We call the $\mathbb{R}$-linear space $\mathbb{R}^n$ \emph{an
algebra}; denoted by $\mathbb A$  if it is endowed with a bilinear
product $\mathbb{A}\times \mathbb{A} \rightarrow \mathbb{A}$
denoted by $(u,v)\mapsto uv$, which is associative and
commutative; that is $u(vw) = (uv)w$ and $ uv = vu$ for all
$u,v,w\in \mathbb{A}$; furthermore, there exists a unit $e \in
\mathbb{A}$, which satisfies $eu=u$ for all $u\in \mathbb{A}$, see
\cite{Pie}.

An element $u\in\mathbb A$ is called \emph{regular} if there
exists $u^{-1}\in\mathbb A$ called \emph{the inverse} of $u$ such
that $u^{-1} u=e$. We also use the notation $e/u$ for $u^{-1}$,
where $e$ is the unit of $\mathbb A$. If $u\in\mathbb A$ is not
regular, then $u$ is called \emph{singular}. $\mathbb A^*$ denotes
the set of all the regular elements of $\mathbb A$. If
$u,v\in\mathbb A$ and $v$ is regular, the quotient $u/v$ means
$uv^{-1}$.

The $\mathbb A$ product between the elements of the canonical
basis $\{e_1,\cdots,e_n\}$ of $\mathbb R^n$ is given by
$e_ie_j=\sum_{k=1}^n c_{ijk}e_k$ where $c_{ijk}\in\mathbb R$ for
$i,j,k\in\{1,\cdots,n\}$ are called \emph{structure constants} of
$\mathbb A$. The \emph{first fundamental representation} of
$\mathbb A$ is the injective linear homomorphism $R:\mathbb A\to
M(n,\mathbb R)$ defined by $R:e_i\mapsto R_i$, where $R_i$ is the
matrix with $[R_i]_{jk}=c_{ikj}$, for $i=1,\cdots,n$.

\subsection{Definition of $\varphi\mathbb A$-differentiability}

We use notation $u=(u_1,\cdots,u_k)$. The usual differential of a
function $f$ will be denoted by $df$.

Let $\mathbb A$ be the linear space $\mathbb R^n$ endowed with an
algebra product. Consider two differentiable functions in the
usual sense $f,\varphi:U\subset\mathbb R^k\to\mathbb R^n$ defined
in an open set $U$. We say $f$ is
$\varphi\mathbb{A}$-\emph{differentiable} on $U$ if there exists a
function $f'_\varphi:U\subset\mathbb R^k\to\mathbb R^n$, that we
call \emph{$\varphi\mathbb{A}$-derivative}, such that for all
$u\in U$
$$
\lim_{\xi\to 0,\,\xi\in\mathbb
R^k}\frac{f(u+\xi)-f(u)-f'_\varphi(u)d\varphi_u(\xi)}{||\xi||}=0,
$$
where $f'_\varphi(u)d\varphi_u(\xi)$ denotes the $\mathbb
A$-product of $f'_\varphi(u)$ and $d\varphi_u(\xi)$. That is, $f$
is $\varphi\mathbb{A}$-differentiable if
$df_u(\xi)=f'_\varphi(u)d\varphi_u(\xi)$ for all $\xi\in\mathbb
R^k$.

If $k=n$ and $\varphi:\mathbb R^n\to\mathbb R^n$ is the identity
transformation $\varphi(x)=x$, the $\varphi\mathbb
A$-differentiability will be called $\mathbb
A$-\emph{differentiability} and the \emph{$\mathbb A$-derivative}
of $f$ will be denoted by $f'$. This last differentiability is
known as Lorch differentiability, see \cite{Lor}.
Differentiability related to commutative and noncommutative
algebras is considered in \cite{JSC}.

\subsection{Algebrizability of planar vector fields}\label{apvf}

The algebrizability ($\mathbb A$-differentiability for some
algebra) of planar vector fields $F=(u,v)$ is characterized in
\cite{AMEE} and \cite{Dyg}. A vector field $F$ is algebrizable on
an open set $\Omega\subset\mathbb R^2$ if and only if $F$
satisfies at least one of the following homogeneous first order
systems of PDEs
\begin{itemize}
    \item [a)] $u_x+\beta v_x-v_y=0$, $u_y-\alpha v_x=0$,
    \item [b)] $u_x+\gamma u_y-v_y=0$, $v_x-\delta u_y=0$, and
    \item [c)] $u_y=0$, $v_x=0$.
\end{itemize}
For case a) we take $\mathbb A=\mathbb A^{2}_1(\alpha,\beta)$, for
b) $\mathbb A=\mathbb A^{2}_2(\gamma,\delta)$, and for c) $\mathbb
A=\mathbb A^{2}_{1,2}$. These systems are called
\emph{Cauchy-Riemann equations} associated with $\mathbb A$
($\mathbb A$-CREs), where $\alpha,\beta,\gamma,\delta\in\mathbb R$
are parameters, see \cite{War2}. For $\mathbb
A^{2}_1(\alpha,\beta)$ the product is
\begin{equation}\label{pa1ab}
    \begin{tabular}{c|c c}
$\cdot$ & $ e_1$& $e_2$ \\ \hline $e_1$&$ e_1$ & $e_2$\\
 $e_2 $&$ e_2 $ & $\alpha e_1+\beta e_2$ \\
\end{tabular}
\end{equation}
hence the unit is $e=e_1$. The structure constants are
\begin{equation}\label{ca1ab}
\begin{tabular}{cccc}
 $c_{111}=1$, & $c_{112}=0$, & $c_{121}=0$, & $c_{122}=1$, \\
 $c_{211}=0$, &  $c_{212}=1$, &  $c_{221}=\alpha$, & $c_{222}=\beta$, \\
\end{tabular}%
\end{equation}
or equivalently, its first fundamental representation is
\begin{equation}\label{ffra1ab}
R(e_1)=\left(%
\begin{array}{cc}
  1 & 0 \\
  0 & 1 \\
\end{array}%
\right),\qquad R(e_2)=\left(%
\begin{array}{cc}
  0 & \alpha \\
  1 & \beta \\
\end{array}%
\right).
\end{equation}
For $\mathbb A^{2}_2(\gamma,\delta)$ the product is
\begin{equation}\label{pa2ab}
    \begin{tabular}{c|c c}
$\cdot$ & $ e_1$& $e_2$ \\ \hline $e_1$ & $\gamma e_1+\delta e_2$  & $e_1$\\
 $e_2 $&$ e_1 $ & $e_2$ \\
\end{tabular}
\end{equation}
hence the unit is $e=e_2$. The structure constants are
\begin{equation}\label{ca2ab}
\begin{tabular}{cccc}
 $c_{111}=\gamma$, & $c_{112}=\delta$, & $c_{121}=1$, & $c_{122}=0$, \\
 $c_{211}=1$, &  $c_{212}=0$, &  $c_{221}=0$, & $c_{222}=1$, \\
\end{tabular}%
\end{equation}
or equivalently, its first fundamental representation is
\begin{equation}\label{ffra2ab}
R(e_1)=\left(%
\begin{array}{cc}
  \gamma & 1 \\
  \delta & 0 \\
\end{array}%
\right),\qquad R(e_2)=\left(%
\begin{array}{cc}
  1 & 0 \\
  0 & 1 \\
\end{array}%
\right).
\end{equation}
For $\mathbb A^{2}_{1,2}$ the product is
\begin{equation}\label{pa3}
    \begin{tabular}{c|c c} $\cdot$ & $ e_1$& $e_2$ \\ \hline $e_1$&$ e_1$ & $0$\\
 $e_2$ & $0$ & $e_2$ \\
\end{tabular}
\end{equation}
hence the unit is $e=e_1+e_2$. The structure constants are
\begin{equation}\label{ca3}
    \begin{tabular}{cccc}
 $c_{111}=1$, & $c_{112}=0$, & $c_{121}=0$, & $c_{122}=0$, \\
 $c_{211}=0$, &  $c_{212}=0$, &  $c_{221}=0$, & $c_{222}=1$, \\
\end{tabular}%
\end{equation}
or equivalently, its first fundamental representation is
\begin{equation}\label{ffra3}
    R(e_1)=\left(%
\begin{array}{cc}
  1 & 0 \\
  0 & 0 \\
\end{array}%
\right),\qquad R(e_2)=\left(%
\begin{array}{cc}
  0 & 0 \\
  0 & 1 \\
\end{array}%
\right).
\end{equation}

\subsection{On $\varphi\mathbb{A}$-differenciability}\label{varphi:derivative}

The $\varphi\mathbb{A}$-derivative $f'_\varphi(u)$ is unique if
$d\varphi_u(\mathbb R^k)$ is not contained in the singular set of
$\mathbb A$. The function $\varphi:U\subset\mathbb R^k\to\mathbb
R^n$ is $\varphi\mathbb A$-differentiable and
$\varphi'_\varphi(u)=e$ for all $u\in U$ where $e\in\mathbb A$ is
the unit. Also, the $\mathbb A$-combinations (linear) and $\mathbb
A$-products of $\varphi\mathbb{A}$-differentiable functions are
$\varphi\mathbb{A}$-differentiable functions and they satisfy the
usual rules of differentiation. In the same way if $f$ is
$\varphi\mathbb{A}$-differentiable and has image in the regular
set of $\mathbb A$, then the function $\frac{e}{f^n}$ is
$\varphi\mathbb{A}$-differentiable for $n\in\{1,2,\cdots\}$, and
\begin{equation}\label{dnsf}
\left(\frac{e}{f^n}\right)'_{\varphi}=\frac{-nf'_\varphi}{f^{n+1}}.
\end{equation}

A $\varphi\mathbb{A}$-\emph{polynomial function} $p:\mathbb
R^k\to\mathbb R^n$ is defined by
\begin{equation}\label{Apolynomial:function}
    p(u)=c_0+c_1\varphi(u)+c_2(\varphi(u))^{2}+\cdots+c_m(\varphi(u))^m
\end{equation}
where $c_0,c_1,\cdots,c_m\in\mathbb A$ are constants and the
variable $u$ represent the variable in $\mathbb R^k$. A
$\varphi\mathbb{A}$-\emph{rational function} is defined by a
quotient of two $\varphi\mathbb{A}$-polynomial functions. Then,
$\varphi\mathbb{A}$-\emph{polynomial functions} and
$\varphi\mathbb{A}$-\emph{rational functions} are $\varphi\mathbb
A$-differentiable and the usual rules of differentiation are
satisfied for the $\varphi\mathbb{A}$-derivative.

In general, the rule of chain does not have sense since $\varphi$
is $\varphi\mathbb{A}$-differentiable however the composition
$\varphi\circ\varphi$ only is defined when $k=n$. Even in this
case the rule of the chain can not be verified. Suppose that
$\varphi$ is a linear isomorphism and that the rule of chain is
satisfied. Thus the Jacobian matrix of $\varphi\circ\varphi$
satisfies
$$
J(\varphi\circ\varphi)=MM=R((\varphi\circ\varphi)'_{\varphi})M,
$$
where $M$ is the matrix associated with $\varphi$ respect to the
canonical basis of $\mathbb R^n$ and $R$ is the first fundamental
representation of $\mathbb A$. Then
\begin{equation}\label{jm}
J(\varphi\circ\varphi)M^{-1}=M\in R(\mathbb A).
\end{equation}
Therefore, if $M$ has determinant $\det(M)\neq 0$ and $M\notin
R(\mathbb A)$, the rule of chain is not valid for the
$\varphi\mathbb{A}$-differentiability. By (\ref{jm}) $M\in
R(\mathbb A)$.

We have the following first version of the rule of chain.

\blem\label{l1} If $g:\Omega\subset\mathbb R^n\to\mathbb R^n$ is
$\mathbb A$-differentiable with $\mathbb A$-derivative $g'$,
$f:U\subset\mathbb R^k\to\mathbb R^n$ is
$\varphi\mathbb{A}$-differentiable, and $f(U)\subset\Omega$, then
$g\circ f:U\subset\mathbb R^k\to\mathbb R^n$ is a
$\varphi\mathbb{A}$-differentiable function with
$\varphi\mathbb{A}$-derivative
$$
(g\circ f)'_\varphi=(g'\circ f) f'_\varphi.
$$
\elem \noindent\textbf{Proof.} The function $g\circ f$ is
differentiable in the usual sense and
$$
d(g\circ
f)_u(\xi)=dg_{f(u)}df_u(\xi)=g'(f(u))f'_\varphi(u)d\varphi_u(\xi).\qquad\Box
$$

Lemma \ref{l1} has the following converse: each $\varphi\mathbb
A$-differentiable function $f$ can be expressed as $g\circ\varphi$
where $g$ is an $\mathbb A$-algebrizable vector field, as we can
see in the following lemma.

\blem\label{l2} If $\varphi:U\subset\mathbb R^n\to\mathbb R^n$ is
a diffeomorphism defined on an open set $U$ and $f:U\subset\mathbb
R^n\to\mathbb R^n$ a is $\varphi\mathbb A$-differentiable on $U$,
then there exists an $\mathbb A$-differentiable vector field $g$
such that $f(u)=g\circ\varphi(u)$ for all $u\in U$, and
$g'({\varphi(u)})=f'_\varphi(u)$. \elem \noindent\textbf{Proof.}
Define $g=f\circ\varphi^{-1}$, thus
$$
dg_{\varphi(u)}=df_{p}d\varphi^{-1}_{\varphi(u)}=f'_\varphi(u)d\varphi_ud\varphi^{-1}_{\varphi(u)}=f'_\varphi(u).
$$
This means that $g$ is $\mathbb A$-differentiable at $\varphi(u)$
and its $\mathbb A$-derivative is $g'(\varphi(u))=f'_\varphi(u)$.
$\Box$

We have the following proposition.

\bprop\label{pl2} Let $\varphi:U\subset\mathbb R^n\to\mathbb R^n$
be a diffeomorphism defined on an open set $U$. The following
three statements are equivalent
\begin{enumerate}
    \item [a)] $f:U\subset\mathbb R^n\to\mathbb R^n$ is $\varphi\mathbb A$-differentiable on $U$.
    \item [b)] $g=f\circ\varphi^{-1}$ is $\mathbb{A}$-differentiable.
    \item [c)] $f$ is differentiable in the usual sense on $U$ and $Jf_u (J\varphi_{u})^{-1}\in R(\mathbb{A})$
    for all $u\in U$.
\end{enumerate}
\eprop \noindent\textbf{Proof.} Suppose a), by Lemma \ref{l2} we
have b).

Suppose b), then $g=f\circ\varphi^{-1}$ is
$\mathbb{A}$-differentiable. Since $\varphi$ is a diffeomorphism
$f=g\circ\varphi$ is differentiable in the usual sense, the rule
of the chain gives $dg_{\varphi(u)}=df_u\circ
d\varphi^{-1}_{\varphi(u)}$, and $Jg_{\varphi(u)}=Jf_u
J\varphi^{-1}_{\varphi(u)}\in
 R(\mathbb{A})$. That is, b) implies c).

Suppose c). Since $f$ is differentiable in the usual sense
$Jg_{\varphi(u)}=Jf_u J\varphi^{-1}_{\varphi(u)}\in R(\mathbb{A})$
implies $Jf_u=Jg_{\varphi(u)}J\varphi_u$. That is,
$df_u=dg_{\varphi(u)}d\varphi_u$. Thus, $f$ is
$\varphi\mathbb{A}$-differentiable. $\Box$

\bcor\label{cvphiapvf} Let $f(x,y)=(u(x,y),v(x,y))$ be a vector
field for which there exists a diffeomorphism
$\phi(s,t)=(x(s,t),y(s,t))$ that is $\varphi=\phi^{-1}$ and
suppose that some of the following conditions are satisfied:
\begin{itemize}
    \item [a)] There exist constants $\alpha$ and $\beta$ such that
    $$\begin{array}{c}
            u_xx_s+u_yy_s+\beta(v_xx_s+v_yy_s)-(v_xx_t+v_yy_t)=0,\\
            u_xx_t+u_yy_t-\alpha(v_xx_s+v_yy_s)=0. \\
          \end{array}
          $$
    \item [b)] There exist constants $\gamma$ and $\delta$ such that
    $$\begin{array}{c}
            u_xx_s+u_yy_s+\gamma(u_xx_t+u_yy_t)-(v_xx_t+v_yy_t)=0,\\
            v_xx_s+v_yy_s-\delta(u_xx_t+u_yy_t)=0. \\
          \end{array}
          $$
    \item [c)] $u_xx_t+u_yy_t=0$ and $v_xx_s+v_yy_s=0$.
\end{itemize}
\ecor In case a) we take $\mathbb{A}=\mathbb
A^{2}_1(\alpha,\beta)$, in b) $\mathbb{A}=\mathbb
A^{2}_2(\gamma,\delta)$, and in c) $\mathbb{A}=\mathbb
A^{2}_{1,2}$. Then $f$ is
$\varphi\mathbb{A}$-differentiable.\\
\noindent\textbf{Proof.} In the three cases the systems of partial
differential equations are the generalized Cauchy-Riemann
equations given in Section \ref{apvf} for $g=f\circ\varphi^{-1}$,
then $g$ is $\mathbb{A}$-differentiable. Thus, by Proposition
\ref{pl2} $f$ is $\varphi\mathbb{A}$-differentiable. $\Box$

We also have the following second version of the rule of chain.

\blem\label{l3} If $\varphi:U\subset\mathbb R^k\to\mathbb R^n$ is
differentiable on an open set $U$, $g:V\subset\mathbb
R^l\to\mathbb R^k$ is differentiable on an open set $V$ with
$g(V)\subset U$, and $f:U\subset\mathbb R^k\to\mathbb R^n$ is
$\varphi\mathbb A$-differentiable on $U$, then $h=f\circ g$ is
$\phi\mathbb A$-differentiable on $V$ for $\phi=\varphi\circ g$,
and $h'_\phi(v)=f'_{\varphi}(g(v))$ \elem \noindent\textbf{Proof.}
We have
$$
dh_v=d(f\circ g)_v=df_{g(v)}dg_v=f'_\varphi(g(v))
d\varphi_{g(v)}dg_v=f'_{\varphi}(g(v))d\phi_v.
$$
Thus, $h'_\phi(v)=f'_{\varphi}(g(v))$. $\Box$

\subsection{On $\varphi\mathbb{A}$-differenciability of quadratic planar vector fields}

All planar vector fields $f$ having the form
$f(x,y)=(a_0+a_1x+a_2y,b_0+b_1x+b_2y)$ are $\varphi\mathbb
A$-differentiable for $\varphi(x,y)=(a_1x+a_2y,b_1x+b_2y)$ since
in this case $f(x,y)=(a_0,b_0)+\varphi(x,y)$, and the identity
function is $\mathbb A$-differentiable with respect to each
two-dimensional algebra $\mathbb A$.

We consider the case of quadratic planar vector fields
\begin{equation}\label{qpvf}
    f(x,y)=(a_0+a_1x+a_2y+a_3x^{2}+a_4xy+a_5y^{2},b_0+b_1x+b_2y+b_3x^{2}+b_4xy+b_5y^{2}).
\end{equation}
In the following propositions we give conditions in order to
determine the $\varphi\mathbb A$-differentiability of the vector
field (\ref{qpvf}) with respect to some linear function
$\varphi:\mathbb R^{2}\to\mathbb R^{2}$ and an algebra $\mathbb
A=\mathbb A^{2}_1(\alpha,\beta)$, $\mathbb A=\mathbb
A^{2}_2(\gamma,\delta)$, and $\mathbb A=\mathbb A_{1,2}$,
respectively.

\bprop\label{algebrizabilidadplanar:caso1} Consider $f$ given in
(\ref{qpvf}) and the matrix
\begin{equation}\label{mat6p4c1}
    M_6(b_i,\alpha,\beta)=
    \left(%
\begin{array}{cccc}
\beta b_1+a_1 & -b_1 & \beta b_2+a_2 & -b_2 \\
 2\beta b_3+2a_3 & -2b_3 & \beta b_4+a_4 & -b_4 \\
  \beta b_4+a_4 & -b_4 & 2\beta b_5+2a_5 & -2b_5 \\
  \alpha b_1 & -a_1 & \alpha b_2 & -a_2 \\
  2\alpha b_3 & -2a_3 & \alpha b_4 & -a_4 \\
  \alpha b_4 & -a_4 & 2\alpha b_5 & -2a_5 \\
\end{array}%
\right).
\end{equation}
If Range $R(M_6(b_i,\alpha,\beta))$ of $M_6(b_i,\alpha,\beta)$
satisfies $R(M_6(b_i,\alpha,\beta))\leq 3$, $ad-bc\neq 0$, and
$v=(a,b,c,d)$ is a vector in the orthogonal complement of the
linear space spanned by the vectors defined by the files of
$M_6(b_i,\alpha,\beta)$, then $f$ is $\varphi\mathbb
A$-differentiable for $\varphi(s,t)=\frac{1}{ad-bc}(ds-bt,at-cs)$
and $\mathbb A=\mathbb A^{2}_1(\alpha,\beta)$.\eprop
\noindent\textbf{Proof.} The Jacobian matrix of $f$ is given by
$$
Jf=\left(%
\begin{array}{cc}
  a_1 & a_2 \\
  b_1 & b_2 \\
\end{array}%
\right)+\left(%
\begin{array}{cc}
  2a_3 & a_4 \\
  2b_3 & b_4 \\
\end{array}%
\right)x+\left(%
\begin{array}{cc}
  a_4 & 2a_5 \\
  b_4 & 2b_5 \\
\end{array}%
\right)y,
$$
and
\begin{eqnarray*}
  Jf\left(%
\begin{array}{cc}
  a & b \\
  c & d \\
\end{array}%
\right) &=& \left(%
\begin{array}{cc}
  a_1a+a_2c & a_1b+a_2d \\
  b_1a+b_2c  & b_1b+b_2d \\
\end{array}%
\right)+\left(%
\begin{array}{cc}
  2a_3a+a_4c & 2a_3b+a_4d \\
  2b_3a+b_4c & 2b_3b+b_4d \\
\end{array}%
\right)x \\
    &+& \left(%
\begin{array}{cc}
  a_4a+2a_5c & a_4b+2a_5d \\
  b_4a+2b_5c & b_4b+2b_5d \\
\end{array}%
\right)y.
\end{eqnarray*}
Since $v=(a,b,c,d)$, $ad-bc\neq 0$, is in the orthogonal
complement of the linear space spanned by the vectors defined by
the files of $M_6(b_i,\alpha,\beta)$, the inner products of these
vectors by $v$ give the system
\begin{eqnarray*}
 0 =(\beta b_1+a_1)a  -b_1b + (\beta b_2+a_2)c -b_2d &=&a_1a+a_2c+\beta(b_1a+b_2c)-(b_1b+b_2d) ,\\
 0 =(2\beta b_3+2a_3)a -2b_3b+ (\beta b_4+a_4)c -b_4d &=& 2a_3a+a_4c+\beta(2b_3a+b_4c)-(2b_3b+b_4d),\\
 0 = (\beta b_4+a_4)a -b_4b+ (2\beta b_5+2a_5)c -2b_5d &=& a_4a+2a_5c+\beta(b_4a+2b_5c)-(b_4b+2b_5d),\\
 0 = \alpha b_1a -a_1b+ \alpha b_2c -a_2d &=& a_1b+a_2d-\alpha(b_1a+b_2c),\\
 0 = 2\alpha b_3a -2a_3b+ \alpha b_4c -a_4d &=& 2a_3b+a_4d -\alpha(2b_3a+b_4c),\\
 0 = \alpha b_4a -a_4b+ 2\alpha b_5c -2a_5d &=& a_4b+2a_5d-\alpha(b_4a+2b_5c).
\end{eqnarray*}
These equations imply equations given in a) of Corollary
\ref{cvphiapvf}, then $f$ is $\varphi\mathbb A$-differentiable for
$\mathbb A=\mathbb A^{2}_1(\alpha,\beta)$. $\Box$

Most practical conditions for determining the $\varphi\mathbb
A$-differentiability of a quadratic planar vector field for some
algebra $\mathbb A$ of the type $\mathbb A^{2}_1(\alpha,\beta)$,
are given in the following corollary.

\bcor\label{ciaab} Consider $f$ given in (\ref{qpvf}) and the
matrix
\begin{equation}\label{mat4p4c1}
    M_4(b_i,\alpha,\beta)=
    \left(%
\begin{array}{cccc}
 2\beta b_3+2a_3 & -2b_3 & \beta b_4+a_4 & -b_4 \\
  \beta b_4+a_4 & -b_4 & 2\beta b_5+2a_5 & -2b_5 \\
  2\alpha b_3 & -2a_3 & \alpha b_4 & -a_4 \\
  \alpha b_4 & -a_4 & 2\alpha b_5 & -2a_5 \\
\end{array}%
\right).
\end{equation}
Let $\alpha$ and $\beta$ be parameters such that $\det
M_4(b_i,\alpha,\beta)= 0$, and $v=(a,b,c,d)$ a vector in the
orthogonal complement of the linear space spanned by the  vectors
defined by the files of $M_4(b_i,\alpha,\beta)$, and $ad-bc\neq
0$, then the planar homogeneous vector
\begin{equation}\label{phvf}
    g(x,y)=(a_3x^{2}+a_4xy+a_5y^{2},b_3x^{2}+b_4xy+b_5y^{2})
\end{equation}
is $\varphi\mathbb A$-differentiable with respect to
$\varphi(s,t)=\frac{1}{ad-bc}(ds-bt,at-cs)$ and $\mathbb A=\mathbb
A^{2}_1(\alpha,\beta)$. Suppose in addition that $v$ is in the
orthogonal complement of the linear space spanned by the vectors
defined by the files of matrix
\begin{equation}\label{mat2p4c1}
    M_2(b_i,\alpha,\beta)=
    \left(%
\begin{array}{cccc}
\beta b_1+a_1 & -b_1 & \beta b_2+a_2 & -b_2 \\
  \alpha b_1 & -a_1 & \alpha b_2 & -a_2 \\
\end{array}%
\right).
\end{equation}
Thus, $f$ is $\varphi\mathbb A$-differentiable for
$\varphi(s,t)=\frac{1}{ad-bc}(ds-bt,at-cs)$ and $\mathbb A=\mathbb
A^{2}_1(\alpha,\beta)$. \ecor

\bprop\label{algebrizabilidadplanar:caso2} Consider $f$ given in
(\ref{qpvf}) and the matrix
\begin{equation}\label{mat6p4c2}
    M_6(b_i,\gamma,\delta)=
    \left(%
\begin{array}{cccc}
  a_1 & \gamma a_1-b_1    & a_2 & \gamma a_2-b_2 \\
 2 a_3 & 2\gamma a_3-2b_3 & a_4 &  \gamma a_4-b_4 \\
   a_4 &  \gamma a_4-b_4 & 2a_5 & 2\gamma a_5-2b_5 \\
   b_1 & -\delta a_1      & b_2 & -\delta a_2 \\
  2 b_3 & -2\delta a_3    & b_4 & -\delta a_4 \\
   b_4 & -\delta a_4     & 2 b_5 & -2\delta a_5 \\
\end{array}%
\right).
\end{equation}
If Range $R(M_6(b_i,\gamma,\delta))$ of $M_6(b_i,\gamma,\delta)$
satisfies $R(M_6(b_i,\gamma,\delta))\leq 3$, $ad-bc\neq 0$, and
$v=(a,b,c,d)$ is a vector in the orthogonal complement of the
linear space spanned by the vectors defined by the files of
$M_6(b_i,\gamma,\delta)$, then $f$ is $\varphi\mathbb
A$-differentiable for $\varphi(s,t)=\frac{1}{ad-bc}(ds-bt,at-cs)$
and $\mathbb A=\mathbb A^{2}_2(\gamma,\delta)$.\eprop
\noindent\textbf{Proof.} The proof is similar to that of Theorem
\ref{algebrizabilidadplanar:caso1}.

Most practical conditions for determining the $\varphi\mathbb
A$-differentiability of a quadratic planar vector field for some
algebra $\mathbb A$ of the type $\mathbb A^{2}_2(\gamma,\delta)$,
are given in the following corollary.

\bcor Consider $f$ given in (\ref{qpvf}) and the matrix
\begin{equation}\label{mat4p4c2}
    M_4(b_i,\gamma,\delta)=
    \left(%
\begin{array}{cccc}
 2 a_3 & 2\gamma a_3-2b_3 & a_4 &  \gamma a_4-b_4 \\
   a_4 &  \gamma a_4-b_4 & 2a_5 & 2\gamma a_5-2b_5 \\
  2 b_3 & -2\delta a_3    & b_4 & -\delta a_4 \\
   b_4 & -\delta a_4     & 2 b_5 & -2\delta a_5 \\
\end{array}%
\right).
\end{equation}
Let $\gamma$ and $\delta$ be parameters such that $\det
M_4(b_i,\gamma,\delta)= 0$, and $v=(a,b,c,d)$ a vector in the
orthogonal complement of the linear space spanned by the  vectors
defined by the files of $M_4(b_i,\gamma,\delta)$, and $ad-bc\neq
0$, then the planar homogeneous vector field
\begin{equation}\label{phvf2}
    g(x,y)=(a_3x^{2}+a_4xy+a_5y^{2},b_3x^{2}+b_4xy+b_5y^{2})
\end{equation}
is $\varphi\mathbb A$-differentiable with respect to
$\varphi(s,t)=\frac{1}{ad-bc}(ds-bt,at-cs)$ and $\mathbb A=\mathbb
A^{2}_2(\gamma,\delta)$. Suppose in addition that $v$ is in the
orthogonal complement of the linear space spanned by the vectors
defined by the files of matrix
\begin{equation}\label{mat2p4c2}
    M_2(b_i,\gamma,\delta)=
    \left(%
\begin{array}{cccc}
   a_1 & \gamma a_1-b_1    & a_2 & \gamma a_2-b_2 \\
   b_1 & -\delta a_1       & b_2 & -\delta a_2 \\
\end{array}%
\right).
\end{equation}
Thus, $f$ is $\varphi\mathbb A$-differentiable for
$\varphi(s,t)=\frac{1}{ad-bc}(ds-bt,at-cs)$ and $\mathbb A=\mathbb
A^{2}_2(\gamma,\delta)$. \ecor

\bprop\label{algebrizabilidadplanar:caso3} Consider $f$ given in
(\ref{qpvf}) and the matrix
\begin{equation}\label{mat6p4c3}
    M_6(b_i)=
    \left(%
\begin{array}{cccc}
0 & a_1 & 0 & a_2 \\
 0 & 2a_3 & 0 & a_4 \\
  0 & a_4 & 0 & 2a_5 \\
   b_1 & 0 &  b_2 & 0 \\
  2 b_3 & 0 &  b_4 & 0 \\
   b_4 & 0 & 2 b_5 & 0 \\
\end{array}%
\right).
\end{equation}
If Range $R(M_6(b_i))$ of $M_6(b_i)$ satisfies $R(M_6(b_i))\leq
3$, $ad-bc\neq 0$, and $v=(a,b,c,d)$ is a vector in the orthogonal
complement of the linear space spanned by the vectors defined by
the files of $M_6(b_i)$, then $f$ is $\varphi\mathbb
A$-differentiable for $\varphi(s,t)=\frac{1}{ad-bc}(ds-bt,at-cs)$
and $\mathbb A=\mathbb A^{2}_{1,2}$.\eprop

Most practical conditions for determining the $\varphi\mathbb
A$-differentiability of a quadratic planar vector field for some
algebra $\mathbb A$ of the type $\mathbb A^{2}_{1,2}$, are given
in the following corollary.

\bcor Consider $f$ given in (\ref{qpvf}) and the matrix
\begin{equation}\label{mat4p4c3}
    M_4(b_i)=
    \left(%
\begin{array}{cccc}
 0 & 2a_3 & 0 & a_4 \\
  0 & a_4 & 0 & 2a_5 \\
  2 b_3 & 0 &  b_4 & 0 \\
   b_4 & 0 & 2 b_5 & 0 \\
\end{array}%
\right).
\end{equation}
If $\det M_4(b_i)= 0$ and $v=(a,b,c,d)$ is a vector in the
orthogonal complement of the linear space spanned by the vectors
defined by the files of $M_4(b_i)$, and $ad-bc\neq 0$, then the
planar homogeneous vector field
\begin{equation}\label{phvf3}
    g(x,y)=(a_3x^{2}+a_4xy+a_5y^{2},b_3x^{2}+b_4xy+b_5y^{2})
\end{equation}
is $\varphi\mathbb A$-differentiable with respect to
$\varphi(s,t)=\frac{1}{ad-bc}(ds-bt,at-cs)$ and $\mathbb
A^{2}_{1,2}$. Suppose in addition that $v$ is in the orthogonal
complement of the linear space spanned by the vectors defined by
the files of matrix
\begin{equation}\label{mat2p4c3}
    M_2(b_i)=
    \left(%
\begin{array}{cccc}
   0 & a_1 & 0 & a_2 \\
   b_1 & 0 &  b_2 & 0 \\
\end{array}%
\right).
\end{equation}
Thus, $f$ is $\varphi\mathbb A$-differentiable for
$\varphi(s,t)=\frac{1}{ad-bc}(ds-bt,at-cs)$ and $\mathbb
A^{2}_{1,2}$. \ecor

\subsection{Cauchy-Riemann equations for the $\varphi\mathbb{A}$-differentiability}

The canonical basis of $\mathbb R^k$ and $\mathbb R^n$ will be
denoted by $\{e_1,\cdots,e_k\}$ and $\{e_1,\cdots,e_n\}$,
respectively, according to the context of the uses it will be
determined if $e_i$ belongs to $\mathbb R^k$ or to $\mathbb R^n$.
The directional derivatives of a function $f$ with respect to a
direction $u\in\mathbb R^k$ (a direction $u$ is a vector with
$\|u\|=1$) is denoted by $f_{u}=(f_{1u},\cdots,f_{nu})$, and the
directional derivative of $f$ with respect to $e_i$ by
$$f_{u_i}=f_{1u_i}e_1+\cdots+f_{nu_i}e_n.$$

The \emph{Cauchy-Riemann equations} for $(\varphi,\mathbb A)$
($\varphi\mathbb A$-\emph{CREs}) is the linear system of $n(k-1)!$
PDEs obtained from
\begin{equation}\label{gcre}
      d\varphi(e_j)f_{u_i} =d\varphi(e_i)f_{u_j}
\end{equation}
for $i,j\in\{1,\cdots,k\}$. For $i=1,\cdots,k$ suppose
$\varphi=(\varphi_1,\cdots,\varphi_n)$, then
\begin{equation}\label{phiei}
    d\varphi(e_i)=\varphi_{u_i}=\sum_{l=1}^n\varphi_{lu_i}e_l.
\end{equation}

In the following theorem the $\varphi\mathbb A$-CREs are given.

\bthm Let $f=(f_1,\cdots,f_n)$ be an $\varphi\mathbb
A$-differentiable function. Thus, the $\varphi\mathbb A$-CREs are
given by
\begin{equation}\label{ecrkn}
    \sum_{m=1}^n\sum_{l=1}^n(f_{mu_i}\varphi_{lu_j}-f_{mu_j}\varphi_{lu_i})C_{lmq}=0
\end{equation}
for $1\leq i<j\leq k$ and $q=1,\cdots,n$, which is a system of
$n(k-1)!$ partial differential equations.  \ethm
\noindent\textbf{Proof.} The equalities (\ref{gcre}) and
(\ref{phiei}) give
$$
\sum_{q=1}^n\left(\sum_{m=1}^n\sum_{l=1}^n(f_{mu_i}\varphi_{lu_j}-f_{mu_j}\varphi_{lu_i})C_{lmq}\right)e_q=0.\qquad\Box
$$

The directional derivatives of $\varphi\mathbb{A}$-differentiable
functions are given in the following lemma.

\blem\label{derdiru} If $f$ is $\varphi\mathbb{A}$-differentiable,
for each direction $u\in\mathbb R^k$ we have
\begin{equation}\label{ddfu}
    f_{u}=f'_\varphi d\varphi_u.
\end{equation}
\elem \noindent\textbf{Proof.} The proof is obtained directly from
the $\varphi\mathbb{A}$-differentiability of $f$. $\Box$

The $\varphi\mathbb A$-differentiability implies the
$\varphi\mathbb A$-CREs, as we see in the following proposition.

\bprop\label{diecr} Let $f:U\subset\mathbb R^k\to\mathbb R^n$ be a
differentiable function in the usual sense on an open set $U$, and
$k\in\{2,\cdots,n\}$. Thus, if $f$ is $\varphi\mathbb
A$-differentiable, then
$d\varphi(e_j)f_{u_i}=d\varphi(e_i)f_{u_j}$. That is, the
components of $f$ satisfy the $\varphi\mathbb A$-CREs. \eprop
\noindent\textbf{Proof.} By using (\ref{ddfu}) we have
$f_{u_i}=f'_\varphi d\varphi(e_i)$ and $f_{u_j}=f'_\varphi
d\varphi(e_j)$. Then
$$
d\varphi(e_j) f_{u_i} = d\varphi(e_j) f'_\varphi d\varphi(e_i)=
d\varphi(e_i)f'_\varphi
d\varphi(e_j)=d\varphi(e_i)f_{u_j}.\qquad\Box
$$

We say $\varphi$ \emph{has an $\mathbb A$-regular direction} $\xi$
if $\xi:U\to\mathbb S^1$ is a function $u\mapsto\xi_u$ such that
$d\varphi_u(\xi_u)$ is a regular element of $\mathbb A$ for all
$u\in U$, where $\mathbb S^1\subset\mathbb R^k$ denotes the unit
sphere centered at the origin. If $\varphi$ has an $\mathbb
A$-regular direction, Proposition \ref{diecr} has a converse, as
we can see in the following theorem.

\bthm\label{tecridu} Let $f:U\subset\mathbb R^k\to\mathbb R^n$ be
a differentiable function in the usual sense on an open set $U$,
and $k\in\{2,\cdots,n\}$. Suppose that $\varphi$ has regular
directions on $U$. Thus, if the components of $f$ satisfies the
$\varphi\mathbb A$-CREs, then $f$ is $\varphi\mathbb
A$-differentiable. \ethm \noindent\textbf{Proof.} Let $\xi$ be a
regular direction of $\varphi$. Since the components of $f$
satisfies the $\varphi\mathbb A$-CREs we have that
$d\varphi(e_j)f_{u_i}=d\varphi(e_i)f_{u_j}$ for $1\leq i,j\leq k$.
Thus,
\begin{eqnarray*}
  d\varphi(u)f_{u_i} &=& f_{u_i}\sum_{j=1}^nu_jd\varphi(e_j)=\sum_{j=1}^nu_jd\varphi(e_j)f_{u_i} \\
   &=& \sum_{j=1}^nu_jd\varphi(e_i)f_{u_j}=\sum_{j=1}^nu_jf_{u_j}d\varphi(e_i)\\
   &=& f_u d\varphi(e_i).
\end{eqnarray*}
Then, $f_{u_i}=\frac{f_\xi}{d\varphi(\xi)}d\varphi(e_i)$. We take
$g_\varphi=\frac{f_\xi}{d\varphi(\xi)}$. By proof of Theorem
\ref{diecr} we have that $df(x) =\sum_{i=1}^kx_i f_{u_i}$ and
$g_\varphi d\varphi(x)=\sum_{i=1}^k x_i g_\varphi d\varphi(e_i)$.
Under these conditions we have that $df(x)=g_\varphi d\varphi(x)$
for all $x\in\mathbb R^k$. That is, $f$ is $\varphi\mathbb
A$-differentiable and $f_\varphi'=g_\varphi$. $\Box$

\subsection{Examples for the complex field}

In the following examples $\mathbb A$ is the complex filed
$\mathbb C$.

We first consider $\varphi(x,y)=(y,x)$.

\bejem Let $\varphi:\mathbb R^2\to\mathbb R^2$ be given by
$\varphi(x,y)=(y,x)$. The CREs are given by
$$
\varphi(e_2)(u,v)_x=\varphi(e_1)(u,v)_y.
$$
Then
$$
e_1(u_x,v_x)=e_2(u_y,v_y)=(-v_y,u_y),
$$
from which we obtain the $\varphi\mathbb A$-CREs for the
$\varphi\mathbb A$-differentiability
$$
u_x=-v_y,\qquad v_x=u_y.
$$

The function $f(x,y)=(y^2-x^2,2xy)$ satisfies
$f(x,y)=(\varphi(x,y))^2$. In this case we have $u(x,y)=y^2-x^2$
and $v(x,y)=2xy$, and they satisfy the $\varphi\mathbb A$-CREs.
\eejem

Now, we consider $\varphi(x,y)=(y,0)$.

\bejem Let $\varphi:\mathbb R^2\to\mathbb R^2$ be given by
$\varphi(x,y)=(y,0)$. The CREs are given by
$$
\varphi(e_2)(u,v)_x=\varphi(e_1)(u,v)_y.
$$
Then $e_1(u_x,v_x)=(0,0)$, from which we obtain the
$\varphi\mathbb A$-CREs for the $\varphi\mathbb
A$-differentiability
$$
u_x=0,\qquad v_x=0.
$$
Thus, the $\varphi\mathbb A$-differentiable functions $f(x,y)$ are
differentiable functions depending only on $y$. \eejem

Now, we consider $\varphi(x,y)=(y,x+y)$.

\bejem Let $\varphi:\mathbb R^2\to\mathbb R^2$ be given by
$\varphi(x,y)=(y,x+y)$. The CREs are given by
$$
\varphi(e_2)(u,v)_x=\varphi(e_1)(u,v)_y.
$$
Then
$$
(e_1+e_2)(u_x,v_x)=e_2(u_y,v_y),
$$
from which we obtain the $\varphi\mathbb A$-CREs for the
$\varphi\mathbb A$-differentiability
$$
v_y=v_x-u_x,\qquad u_y=u_x+v_x.
$$

The function $f(x,y)=(-x^2-2xy,2xy+2y^2)$ satisfies
$f(x,y)=(\varphi(x,y))^2$. In this case we have $u(x,y)=-x^2-2xy$,
$v(x,y)=2xy+2y^2$ and they satisfy the $\varphi\mathbb A$-CREs.
\eejem

\bejem Let $\varphi:\mathbb R^3\to\mathbb R^2$ be defined by
$\varphi(x,y,z)=(x+z,y)$. Thus, the $\varphi\mathbb A$-CREs are
given by
\begin{equation}\label{xzy}
  \begin{array}{cc}
    u_x=v_y, &  v_x=-u_y, \\
    u_x=u_z, &  v_x=v_z.\\
  \end{array}
\end{equation}

Since the components of the function $f(x,y,z)=\varphi(x,y,z)$
satisfy equations (\ref{xzy}) and $\varphi(e_1)=e$, by Theorem
\ref{tecridu} we have that $f$ is $\varphi\mathbb
A$-differentiable. Thus, polynomial functions $p$ of the form
$$
p(x,y,z)\mapsto c_0+c_1\varphi(x,y,z)+\cdots+c_m(\varphi(x,y,z))^m
$$
where $c_0,c_1,\cdots,c_m\in\mathbb A$ are $\varphi\mathbb
A$-differentiable. Moreover, the rational functions obtained by
quotients of these polynomials are $\varphi\mathbb
A$-differentiable on their domain of definition.

The $f_1(x,y,z)=x^2+z^2+2xz-y^2$ and $f_2(x,y,z)=2xy+2yz$ satisfy
the $\varphi\mathbb A$-CREs, thus $f(x,y,z)=(x^2-2xy,2xy+2yz)$ is
$\varphi\mathbb A$-differentiable. Since $f(x,0,0)=x^2(1,0)$, we
have $f(x,y,z)=(\varphi(x,y,z))^2$. \eejem

We consider an example in which $\varphi$ is non linear.

\bejem Let $\varphi:\mathbb R^3\to\mathbb R^2$ be defined by
$\varphi(x,y,z)=(x^2+z,1/y)$. Thus, the $\varphi\mathbb A$-CREs
are given by
\begin{equation}\label{xzy}
  \begin{array}{cc}
    -\frac{1}{y^2}u_x=2x v_y, &  \frac{1}{y^2}v_x=2xu_y, \\
    u_x=2xu_z, &  v_x=2xv_z.\\
  \end{array}
\end{equation}

Since the components of the function $f(x,y,z)=\varphi(x,y,z)$
satisfy equations (\ref{xzy}) and $d\varphi(e_3)=e$, by Theorem
\ref{tecridu} we have that $f$ is $\varphi\mathbb
A$-differentiable. Thus, polynomial functions $p$ of the form
$$
p(x,y,z)\mapsto c_0+c_1\varphi(x,y,z)+\cdots+c_m\varphi^m(x,y,z)
$$
where $c_0,c_1,\cdots,c_m\in\mathbb A$ are $\varphi\mathbb
A$-differentiable. Moreover, the rational functions obtained by
quotients of these polynomials are $\varphi\mathbb
A$-differentiable on their domain of definition. \eejem

The systems of $\varphi\mathbb A$-CREs presented in this section
do not satisfy the conditions given in \cite{War2}, then the
$\varphi\mathbb A$-differentiable functions related to these
examples are not $\mathbb A$-differentiable for all the algebras
$\mathbb A$.

\subsection{Examples for $\mathbb A=\mathbb
A^3_1(p_1,\cdots,p_6)$}\label{sap16}

Consider the algebra $\mathbb A=\mathbb A^3_1(p_1,\cdots,p_6)$
with unit $e=(1,0,0)$ given in Theorem 1.3 of \cite{Dyg3}. That
is, $\mathbb A$ is the linear space $\mathbb{R}^{3}$ endowed with
the product
\begin{equation}\label{algebra1}
  \begin{tabular}{c|ccc}
  $\cdot$ & $e_1$ & $e_2$ & $e_3$ \\
  \hline
  $e_1$ & $e_1$ & $e_2$ & $e_3$ \\
  $e_2$ & $e_2$ & $p_7e_1+p_1e_2+p_2e_3$ & $p_8e_1+p_3e_2+p_4e_3$ \\
  $e_3$ & $e_3$ & $p_8e_1+p_3e_2+p_4e_3$ & $p_9e_1+p_5e_2+p_6e_3$ \\
\end{tabular}
\end{equation}
where $e=e_1$ and $p_1,p_2,\cdots,p_9$ are parameters satisfying
the equalities
\begin{equation}\label{ceq}
    \begin{array}{c}
      p_7 = p_2p_3+p_4^2-p_1p_4-p_2p_6, \\
      p_8= p_2p_5-p_3p_4,\qquad\qquad\,\,\,\,\,\,\\
      p_9= p_3^2+p_4p_5-p_1p_5-p_3p_6, \\
    \end{array}
\end{equation}

In the following examples $\mathbb A$ is the algebra $\mathbb
A^3_1(p_1,\cdots,p_6)$.

\bejem Let $\varphi:\mathbb R^2\to\mathbb R^3$ be defined by
$\varphi(x,y)=(x,y,0)$. Thus, the $\varphi\mathbb A$-CREs are
given by
\begin{equation}\label{xy0}
   \begin{array}{ccc}
  u_y & = & p_7v_x+p_8w_x \\
  v_y & =  & u_x+p_1v_x+p_3w_x \\
  w_y  & = & p_2v_x+p_4w_x \\
\end{array}.
\end{equation}

Since the components of the function $f(x,y)=\varphi(x,y)$ satisfy
equations (\ref{x0y}) and $\varphi(1,0)=e$, by Theorem
\ref{tecridu} we have that $f$ is $\varphi\mathbb
A$-differentiable. Thus, polynomial functions $p$ of the form
$$
p(x,y)\mapsto c_0+c_1\varphi(x,y)+\cdots+c_m(\varphi(x,y))^m
$$
where $c_0,c_1,\cdots,c_m\in\mathbb A$ are $\varphi\mathbb
A$-differentiable. Moreover, the rational functions obtained by
quotients of these polynomials are $\varphi\mathbb
A$-differentiable on their domain of definition. \eejem

\bejem Let $\varphi:\mathbb R^2\to\mathbb R^3$ be defined by
$\varphi(x,y)=(x,0,y)$. Thus, the $\varphi\mathbb A$-CREs are
given by
\begin{equation}\label{x0y}
   \begin{array}{ccc}
  u_y & = & p_8v_x+p_9w_x \\
  v_y & =  & p_3v_x+p_5w_x \\
  w_y  & = & u_x+p_4v_x+p_6w_x \\
\end{array}.
\end{equation}

Since the components of the function $f(x,y)=\varphi(x,y)$ satisfy
equations (\ref{x0y}) and $\varphi(1,0)=e$, by Theorem
\ref{tecridu} we have that $f$ is $\varphi\mathbb
A$-differentiable. Thus, polynomial functions $p$ of the form
$$
p(x,y)\mapsto c_0+c_1\varphi(x,y)+\cdots+c_m(\varphi(x,y))^m
$$
where $c_0,c_1,\cdots,c_m\in\mathbb A$ are $\varphi\mathbb
A$-differentiable. Moreover, the rational functions obtained by
quotients of these polynomials are $\varphi\mathbb
A$-differentiable on their domain of definition. \eejem

\bejem Let $\varphi:\mathbb R^2\to\mathbb R^3$ be defined by
$\varphi(x,y)=(0,x,y)$. Thus, the $\varphi\mathbb A$-CREs are
given by
\begin{equation}\label{0xy}
   \begin{array}{ccc}
  p_8v_x+p_9w_x-p_7v_y-p_8w_y & = & 0 \\
  p_3v_x+p_5w_x-u_y-p_1v_y-p_3w_y & =  & 0 \\
  u_x+p_4v_x+p_6w_x-p_2v_y-p_4w_y & = & 0 \\
\end{array}.
\end{equation}
\eejem

Consider the function $f(x,y)=\varphi(x,y)$, then
$df_{(a,b)}(x,y)=(0,x,y)$. That is, $f$ is $\varphi\mathbb
A$-differentiable and $f'_{\varphi}(x,y)=e$ for all
$(x,y)\in\mathbb R^2$. Thus, polynomial functions $p$ of the form
$$
p(x,y)\mapsto c_0+c_1\varphi(x,y)+\cdots+c_m(\varphi(x,y))^m
$$
where $c_0,c_1,\cdots,c_m\in\mathbb A$ are $\varphi\mathbb
A$-differentiable. Moreover, the rational functions obtained by
quotients of these polynomials are $\varphi\mathbb
A$-differentiable on their domain of definition.

\bejem Any function $f(x,y)=(0,f_2(x,y),f_3(x,y))$ satisfies
$\varphi\mathbb A$-CREs (\ref{0xy}) for $\mathbb A=\mathbb
A^3_1(0,\cdots,0)$ and $\varphi(x,y)=(0,x,y)$, but not all the
differentiable (in the usual sense) functions $f$ having this form
are $\varphi\mathbb A$-differentiable. This can happen because
$\varphi(\mathbb R^2)$ is contained in the singular set of
$\mathbb A$.

Suppose that $f(x,y)=(0,x+y,x-y)$ is $\varphi\mathbb
A$-differentiable. We have that $df_{(0,0)}(x,y)=(0,x+y,x-y)$.
Thus, if $(a,b,c)=f'_{\varphi}(0,0)$, then
$$
(0,x+y,x-y)=(a,b,c)(0,x,y)=(0,ax,ay)
$$
which is a contradiction. Therefore, $f$ satisfies equations
(\ref{0xy}) but $f$ is not $\varphi\mathbb A$-differentiable.
\eejem

The type of $\varphi\mathbb A$-CREs systems considered in this
section are not considered in \cite{War2}.

\section{Linear systems of two first order PDEs}\label{s2}

\subsection{Equivalent systems of two PDEs}

\noindent Consider the linear first order system of PDEs
\begin{equation}\label{eqn:generalcase}
\begin{array}{ccc}
  a_{11}u_x+ a_{12}u_y+a_{13}v_x+a_{14}v_y  & = & f, \\
  a_{21}u_x+ a_{22}u_y+a_{23}v_x+a_{24}v_y & = & g, \\
\end{array}%
\end{equation}
where $a_{ij}$, $f$ and $g$ are differentiable functions of $u$,
$v$, $x$, and $y$. The corresponding homogenous system is given by
\begin{equation}\label{eqn:hcase}
\begin{array}{ccc}
  a_{11}u_x+ a_{12}u_y+a_{13}v_x+a_{14}v_y  & = & 0, \\
  a_{21}u_x+ a_{22}u_y+a_{23}v_x+a_{24}v_y & = & 0. \\
\end{array}%
\end{equation}

We use notation
\begin{equation}\label{eqn:gcmatrix}
   \langle \left(%
\begin{array}{cccc}
  a_{11} & a_{12} & a_{13} & a_{14} \\
 a_{21} & a_{22} & a_{23} & a_{24} \\
\end{array}%
\right):\left(%
\begin{array}{cc}
  u_x & u_y \\
  v_x & v_y \\
\end{array}%
\right)\rangle=\left(%
\begin{array}{c}
  f \\
  g \\
\end{array}%
\right),
\end{equation}
for system (\ref{eqn:generalcase}). The corresponding homogenous
system is expressed by
\begin{equation}\label{eqn:hgcmatrix}
   \langle \left(%
\begin{array}{cccc}
  a_{11} & a_{12} & a_{13} & a_{14} \\
 a_{21} & a_{22} & a_{23} & a_{24} \\
\end{array}%
\right):\left(%
\begin{array}{cc}
  u_x & u_y \\
  v_x & v_y \\
\end{array}%
\right)\rangle=\left(%
\begin{array}{c}
  0 \\
  0 \\
\end{array}%
\right).
\end{equation}
If $A$ is the $2\times 4$ matrix function formed by the $a_{ij}$,
$w=(u,v)$, $F=(f,g)$ and $dw$ denotes the usual derivative, then
(\ref{eqn:gcmatrix}) is written by $\langle A:dw\rangle=F$.

\bdefe If $A_1$, $A_2$ are two $2\times 4$ matrix valued
functions, and $F_1$, $F_2$ are two $2\times 1$ matrix valued
functions, we say that systems  $\langle A_1:dw\rangle=F_1$ and
$\langle A_2:dw\rangle=F_2$ are \emph{equivalent} if there exists
a non singular $2\times 2$ matrix valued function $M=M(x,y)$ such
that $A_2=MB_1$ and $F_2=MF_1$. \edefe

If  $\langle A_1:dw\rangle=F_1$ and  $\langle A_2:dw\rangle=F_2$
are equivalent, then they share the same solutions set.

\bprop Let $\langle A_1:dU\rangle=F_1$ and  $\langle
A_2:dw\rangle=F_2$ be systems that share the same solutions set.
If there exists a $2\times 2$ matrix valued function $M$ such that
$F_1=MF_2$ and there exists a solution $w$ with $dw$ not singular,
then the given systems are equivalent. \eprop
\noindent\textbf{Proof.} Let $w$ be a solution with $dw$ being non
singular, then
$$
\langle A_1-MA_2:dw\rangle=\langle A_1:dw\rangle-\langle
MA_2:dw\rangle=F_1-MF_2=0,
$$
that is, $A_1-MA_2=0$. Thus, the given systems are equivalent.
$\Box$

If $w_p$ is a particular solution of system
(\ref{eqn:generalcase}), then each solution $w$ of system
(\ref{eqn:generalcase}) is expressed by $w=w_h+w_p$, where $w_h$
is a solution of the associated homogenous system
(\ref{eqn:hcase}). Therefore, if the homogenous system results
equivalent to a set of CREs for the $\varphi\mathbb{A}$-derivative
and a particular solution $w_p$ of (\ref{eqn:generalcase}) is
known, then the solutions of (\ref{eqn:generalcase}) have the form
$w=w_h+w_p$, where $w_h$ is a $\varphi\mathbb{A}$-differentiable
function.

\subsection{Homogeneous linear systems of two first order PDEs}

The problem considered in this section is when a given homogenous
linear system of two PDEs is the CREs for the
$\varphi\mathbb{A}$-differentiability for some differentiable
function $\varphi$ and an algebra $\mathbb{A}$. The general system
we consider has the form
\begin{equation}\label{hlspde2}
   \begin{array}{ccc}
  a_{11}u_x+ a_{12}u_y+a_{13}v_x+a_{14}v_y & = & 0\\
 a_{21}u_x+ a_{22}u_y+a_{23}v_x+a_{24}v_y & = & 0 \\
\end{array},
\end{equation}
where $a_{ij}$ are functions of $(x,y)$ for $i=1,2$ and
$j=1,\cdots,4$. Suppose that
$$
(a_{11},a_{12},a_{13},a_{14})\neq
\alpha(a_{21},a_{22},a_{23},a_{24})
$$
for all differentiable scalar functions $\alpha(x,y)$.

We are interested in the existence of an algebra $\mathbb{A}$ and
a scalar function $\varphi(x,y)$ such that system (\ref{hlspde2})
is equivalent to a system which is the CREs for the
$\varphi\mathbb{A}$-derivative.

We have the following propositions and examples.

\bprop\label{tvpasa1} Suppose that $F_1=(-b_{12},b_{11})$ and
$F_2=(-b_{22},b_{21})$ are conservative vector fields with
potential functions $\varphi_i$ for $i=1,2$, such that system
$\langle A:dw\rangle=0$ is equivalent to system
\begin{equation}\label{sal1}
   \begin{array}{ccc}
  b_{11}u_x+ b_{12}u_y+\alpha b_{21}v_x+\alpha b_{22}v_y & = & 0\\
 b_{21}u_x+ b_{22}u_y+(b_{11}+\beta b_{21})v_x+(b_{12}+\beta b_{22})v_y & = & 0 \\
\end{array},
\end{equation}
for some parameters $\alpha,\beta$, then the
$\varphi\mathbb{A}$-differentiable functions are solutions of
$\langle A:dw\rangle=0$, where $\varphi=(\varphi_1,\varphi_2)$,
and $\mathbb{A}=\mathbb{A}^2_1(\alpha,\beta)$ is the algebra given
in Section \ref{s1}. \eprop\noindent\textbf{Proof.} Under the
conditions given above the Jacobian matrix of $\varphi$ is given
by
$$
J\varphi=\left(%
\begin{array}{cc}
  -b_{12} & b_{11} \\
  -b_{22} & b_{21}  \\
\end{array}%
\right).
$$
The $\varphi\mathbb{A}$-CREs are given by
$d\varphi(e_2)(u_x,v_x)=d\varphi(e_1)(u_y,v_y)$, where $d\varphi$
denotes the usual derivative. That is,
$$
(b_{11},b_{21} )(u_x,v_x)=-(b_{12},b_{22})(u_y,v_y).
$$
By using the $\mathbb{A}$-product we obtain system (\ref{sal1}).
$\Box$

\bejem Consider the system
\begin{equation}\label{sal1e}
\begin{array}{ccc}
 \begin{array}{ccc}
    yu_x+xu_y-\alpha xv_x+\alpha yv_y & = & 0\\
 xu_x-yu_y-(y-\beta x)v_x-(x+\beta y)v_y & = & 0 \\
\end{array}.
\end{array}
\end{equation}
Since $F_1=(2x,-2y)$, $F_2=(2y,2x)$ are conservative vector fields
with potential functions $\varphi_1(x,y)=x^{2}-y^{2}$,
$\varphi_2(x,y)=2xy$, by Proposition \ref{diecr} and Theorem
\ref{tecridu} we have for $\varphi=(\varphi_1,\varphi_2)$ and
$\mathbb{A}=\mathbb{A}^{2}_{1}(\alpha,\beta)$ that the set of
$\varphi\mathbb{A}$-differentiable functions is the set of the
solutions of system (\ref{sal1e}).

By Lemma \ref{l2} we have that $\varphi\mathbb{A}$-differentiable
functions $f$ have the form $f=g\circ \varphi$, where $g$ is an
$\mathbb{A}$-differentiable function. So,
$$
f(x,y)=g(x^{2}-y^{2},2xy),
$$
is solution of (\ref{sal1e}) if $g$ is an
$\mathbb{A}$-differentiable function. \eejem

\bprop\label{tvpasa2} Suppose that $F_1=(-b_{14},b_{13})$ and
$F_2=(-b_{24},b_{23})$ are conservative vector fields with
potential functions $\varphi_i$ for $i=1,2$, such that system
$\langle A:dw\rangle=0$ is equivalent to system
\begin{equation}\label{sal2}
   \begin{array}{ccc}
  (\gamma b_{13}+b_{23})u_x+(\gamma b_{14}+b_{24})u_y +b_{13}v_x+b_{14}v_y & = & 0\\
 \delta b_{13}u_x+ \delta b_{14}u_y+b_{23}v_x+b_{24}v_y & = & 0 \\
\end{array},
\end{equation}
for some parameters $\gamma,\delta$, then the
$\varphi\mathbb{A}$-differentiable functions are solutions of
system $\langle A:dw\rangle=0$, where
$\varphi=(\varphi_1,\varphi_2)$ and
$\mathbb{A}=\mathbb{A}^2_2(\gamma,\delta)$ is the algebra given in
Section \ref{s1}. \eprop \noindent\textbf{Proof.} The proof is
similar to that of \ref{tvpasa1}. $\Box$

\bejem Consider the system
\begin{equation}\label{sal2e}
\begin{array}{ccc}
(\gamma y-x)u_x+(\gamma x+y)u_y +yv_x+xv_y & = & 0\\
 \delta yu_x+\delta xu_y-xv_x+yv_y & = & 0 \\
\end{array}.
\end{equation}
Since $F_1=(2x,-2y)$, $F_2=(2y,2x)$ are conservative vector fields
with potential functions $\varphi_1(x,y)=x^{2}-y^{2}$,
$\varphi_2(x,y)=2xy$, by Proposition \ref{diecr} and Theorem
\ref{tecridu} we have for $\varphi=(\varphi_1,\varphi_2)$ and
$\mathbb{A}=\mathbb{A}^{2}_{2}(\gamma,\delta)$ that the set of
$\varphi\mathbb{A}$-differentiable functions is the set of the
solutions of system (\ref{sal2e}).

By Lemma \ref{l2} we have that $\varphi\mathbb{A}$-differentiable
functions $f$ have the form $f=g\circ \varphi$, where $g$ is an
$\mathbb{A}$-differentiable function. So,
$$
f(x,y)=g(x^{2}-y^{2},2xy),
$$
is solution of (\ref{sal2e}) if $g$ is an
$\mathbb{A}$-differentiable function. \eejem

\bprop\label{tvpasa3} Suppose that $F_1=(-b_{2},b_{1})$ and
$F_2=(-b_{4},b_{3})$ are conservative vector fields with potential
functions $\varphi_i$ for $i=1,2$, such that system $\langle
A:dw\rangle=0$ is equivalent to system
\begin{equation}\label{sal3}
   \begin{array}{ccc}
  b_{1}u_x+b_{2}u_y  & = & 0\\
b_{3}v_x+b_{4}v_y & = & 0 \\
\end{array},
\end{equation}
then the $\varphi\mathbb{A}$-differentiable functions are
solutions of system $\langle A:dw\rangle=0$, where
$\mathbb{A}=\mathbb{A}^2_{1,2}$ is the algebra given in Section
\ref{s1}. \eprop \noindent\textbf{Proof.} The proof is similar to
that of \ref{tvpasa1}. $\Box$

\bejem Consider the system
\begin{equation}\label{sdedp3}
\begin{array}{ccc}
  yu_x+xu_y  & = & 0,\\
xv_x-yv_y & = & 0. \\
\end{array}
\end{equation}
Since $F_1=(x,-y)$, $F_2=(y,x)$ are conservative vector fields
with potential functions
$\varphi_1(x,y)=\frac{x^{2}}{2}-\frac{y^{2}}{2}$,
$\varphi_2(x,y)=xy$, by Proposition \ref{diecr} and Theorem
\ref{tecridu} we have for $\varphi=(\varphi_1,\varphi_2)$ and
$\mathbb{A}=\mathbb{A}^{2}_{1,2}$ that the set of
$\varphi\mathbb{A}$-differentiable functions is the set of the
solutions of system (\ref{sdedp3}).

The $\mathbb{A}$-differentiable functions have the form
$g(x,y)=(g_1(x),g_2(y))$ where $g_1$ and $g_2$ are one variable
differentiable functions in the usual sense. By Lemma \ref{l2} we
have that $\varphi\mathbb{A}$-differentiable functions $f$ have
the form $f=g\circ \varphi$, where $g$ is an
$\mathbb{A}$-differentiable function. So,
$$
f(x,y)=\left(g_1\left(\frac{x^{2}}{2}-\frac{y^{2}}{2}\right),g_2(xy)\right),
$$
is solution of (\ref{sdedp3}) if $g_1$ and $g_2$ are
differentiable functions of one variable. \eejem

\section{$\varphi\mathbb A$-differential equations}\label{s3}

\subsection{The Cauchy-integral theorem for the $\varphi\mathbb
A$-differentiability}

If $f:U\subset\mathbb R^k\to\mathbb R^n$ is a $\varphi\mathbb
A$-differentiable function defined in an open set $U$. The
$\varphi\mathbb A$-\emph{line integral} of $f$ is defined by
\begin{equation}\label{pail}
\int_\gamma fd\varphi=\int_\gamma
f(v)d\varphi(v'):=\int_0^{t_1}f(\gamma(s))d\varphi(\gamma'(s))ds,
\end{equation}
where $\gamma$ is a differentiable function of $t$ with values in
$U$, $\gamma(0)=u_0$, $\gamma(t_1)=u$,
$f(\gamma(s))d\varphi(\gamma'(s))$ represents the $\mathbb
A$-product of $f(\gamma(s))$ and
$d\varphi_{\gamma(s)}(\gamma'(s))$, and the right hand of
(\ref{pail}) represents the usual line integral in $\mathbb R^n$.

A version of the Cauchy integral theorem for the $\varphi\mathbb
A$-line integral is given in the following theorem, see \cite{EL1}
and Corollary 10.11 pg. 49 of \cite{JSC} for another version of
the Cauchy-integral theorem relative to algebras.

\bthm\label{ticpd} Let $\varphi:U\subset\mathbb R^k\to\mathbb R^n$
be a $C^2$-function defined on a simply-connected open set $U$ and
$f:U\subset\mathbb R^k\to\mathbb R^n$ a $\varphi\mathbb
A$-differentiable function. If $\gamma$ is a closed differentiable
path contained in $U$, then the $\varphi\mathbb A$-line integral
(\ref{pail}) is equal to zero. \ethm \noindent\textbf{Proof.} We
will show that $fd\varphi(\gamma')=\sum_{q=1}^n\langle
G_q,\gamma'\rangle e_q$, where the $G_q$ are n-dimensional
conservative vector fields. Remember that
$d\varphi(e_j)=\sum_{l=1}^n \varphi_{lu_j}e_l$.

The $\mathbb A$-product of $f$ and $\varphi(\gamma')$ is given by
\begin{eqnarray*}
  fd\varphi(\gamma') &=& \left(\sum_{m=1}^nf_me_m\right)\left(\sum_{j=1}^k\gamma_j'd\varphi(e_j)\right)=
  \sum_{m=1}^n\sum_{j=1}^k\sum_{l=1}^n f_m\gamma_j'\varphi_{lu_j}e_le_m \\
   &=&  \sum_{q=1}^n\left(\sum_{m=1}^n\sum_{j=1}^k\sum_{l=1}^n
   f_m\gamma_j'\varphi_{lu_j}C_{lmq}\right)e_q\\
   &=& \sum_{q=1}^n\left(\langle\sum_{j=1}^k\left(\sum_{m=1}^n\sum_{l=1}^n
   f_m\varphi_{lu_j}C_{lmq}\right)e_j,\sum_{j=1}^k\gamma_j'e_j\rangle\right)e_q,
\end{eqnarray*}
where $\langle\cdot,\cdot\rangle$ denotes the inner product of the
vector field $G_q$ and $\gamma'$, and
$$
G_q=\sum_{j=1}^k\left(\sum_{m=1}^n\sum_{l=1}^n
   f_m\varphi_{lu_j}C_{lmq}\right)e_j
$$
for $q=1,\cdots,n$. By taking the exterior derivative of the dual
1-form of $F_q$, using the $\varphi\mathbb A$-CREs given by
(\ref{ecrkn}), and the commutativity of the second partial
derivatives of the components of $\varphi$, we show that this
1-form is exact. Therefore, $G_q$ is a conservative vector field.
See \cite{EL1} for the proof of this theorem for the $\mathbb
A$-differentiable functions. $\Box$

If $U$ is a simply connected open set containing $u$ and $u_0$,
Theorem \ref{ticpd} permit us to define
$$
\int_{u_0}^u fd\varphi=\int_\gamma fd\varphi,
$$
where $\gamma$ is a differentiable function of $t$ with values in
$U$, $\gamma(0)=u_0$, and $\gamma(t_1)=u$, as in Definition
\ref{pail}.

\bcor\label{ctic} Let $\varphi:U\subset\mathbb R^k\to\mathbb R^n$
be a $C^2$-function on an open set $U$ and $f:U\subset\mathbb
R^k\to\mathbb R^n$ be a $\varphi\mathbb A$-differentiable
function. The vector fields
$$
G_q=\sum_{j=1}^k\left(\sum_{m=1}^n\sum_{l=1}^n
   f_m\varphi_{lu_j}C_{lmq}\right)e_j
$$
for $q=1,\cdots,n$ are conservative, where
$\varphi_{u_j}=\sum_{l=1}^n \varphi_{lu_j}e_l$. \ecor

\bejem Consider the algebra $\mathbb A=\mathbb
A^3_1(-1,\cdots,-1)$ with unit $e=(1,0,0)$ given in Subsection
\ref{sap16} which is given by
\begin{equation}\label{algebra1}
  \begin{tabular}{c|ccc}
  $\cdot$ & $e_1$ & $e_2$ & $e_3$ \\
  \hline
  $e_1$ & $e_1$ & $e_2$ & $e_3$ \\
  $e_2$ & $e_2$ & $e_2+e_3$ & $e_2+e_3$ \\
  $e_3$ & $e_3$ & $e_2+e_3$ & $e_2+e_3$ \\
\end{tabular}.
\end{equation}
The structure constants of $\mathbb A$ are given by
$$
\begin{array}{ccc}
  C_{111}=1, & C_{112}=0, & C_{113}=0, \\
  C_{121}=0, & C_{122}=1, & C_{123}=0,  \\
  C_{131}=0, & C_{132}=0, & C_{133}=1,  \\
   C_{211}=0, & C_{212}=1, & C_{213}=0, \\
  C_{221}=0, & C_{222}=1, & C_{223}=1,  \\
  C_{231}=0, & C_{232}=1, & C_{233}=1,  \\
    C_{311}=0, & C_{312}=0, & C_{313}=1, \\
  C_{321}=0, & C_{322} =1, & C_{323}=1,  \\
  C_{331}=0, & C_{332}=1, & C_{333}=1.  \\
\end{array}%
$$

Let $\varphi(x,y)=(x,y,0)$. The function
$f(x,y,z)=\varphi(x,y)^{-1}$ is $\varphi\mathbb A$-differentiable
and
$$
f(x,y)=\left(\frac{1}{x},\frac{-xy-y^2}{x^3+2x^2y},\frac{y^2}{x^3+2x^2y}\right).
$$

Thus, the conservative vector fields $G_i$ for $i=1,2,3$ are given
by
\begin{eqnarray*}
  G_1 &=&(f_1,0)= \left(\frac{1}{x},0\right),\\
  G_2 &=& (f_2,f_1+f_2+f_3)= \left(\frac{-xy-y^2}{x^3+2x^2y},\frac{x+y}{x^2+2xy}\right), \\
  G_3 &=&
  (f_3,f_2+f_3)=\left(\frac{y^2}{x^3+2x^2y},-\frac{xy}{x^3+2x^2y}\right).
\end{eqnarray*}
\eejem

If $U$ is a simply-connected open set and $f:U\subset\mathbb
R^k\to\mathbb R^n$ is $\varphi\mathbb A$-differentiable on $U$,
then the function
$$
F(u)=\int_{u_0}^u f(v)d\varphi(v'),
$$
for $u_0,u\in U$ is well defined. An \emph{$\varphi\mathbb
A$-antiderivative} of a function $f:U\subset\mathbb R^k\to\mathbb
R^n$ is a function $F:U\subset\mathbb R^k\to\mathbb R^n$ whose
$\varphi\mathbb A$-derivative is given by $F'_\varphi=f$.

For $\varphi\mathbb A$-polynomial functions the $\varphi\mathbb
A$-antiderivative can be calculated in the usual way. The
$\varphi\mathbb A$-line integral of $\varphi\mathbb
A$-differentiable functions gives $\varphi\mathbb
A$-antiderivatives, as we have in the following corollary which is
a generalization of the fundamental theorem of calculus.

\bcor\label{ctfc} Let $\varphi:U\subset\mathbb R^k\to\mathbb R^n$
be a $C^2$-function defined on a simply-connected open set $U$ and
$f:U\subset\mathbb R^k\to\mathbb R^n$ a $\varphi\mathbb
A$-differentiable function. If $u_0,u\in U$ and
$$
F(u)=\int_{u_0}^u f(v)d\varphi(v'),
$$
then $F'_\varphi=f$. \ecor \noindent\textbf{Proof.} We take the
curve $\gamma(t)=u+t\xi$ joining $u$ and $u+\xi$, thus
$\gamma'(t)=\xi$. The rest of the proof is a consequence of
Theorem \ref{ticpd}. $\Box$

\subsection{Existence and uniqueness of solutions}\label{spade}

Let $F:\Omega\subset\mathbb R^n\to\mathbb R^n$ be a vector field
defined on an open set $\Omega$. A $\varphi\mathbb
A$-\emph{differential equation} is
\begin{equation}\label{apode}
    w'_\varphi=F(w),\qquad w(\tau_0)=w_0,
\end{equation}
finding a solution is understand as the problem of finding a
$\varphi\mathbb A$-differentiable function
$w:V_{\tau_0}\subset\mathbb R^k\to\mathbb R^n$ defined in a
neighborhood $V_{\tau_0}$ of $\tau_0$ such that
$dw_\tau=F(w(\tau))d\varphi_\tau$ for all $\tau\in V_{u_0}$, and
satisfying the initial condition $w(\tau_0)=w_0$.

We have the following Existence and Uniqueness Theorem for
$\mathbb A$-algebrizable vector fields and $\varphi\mathbb
A$-differential equations.

\bthm\label{tcvcvaa} Let $\varphi:U\subset\mathbb R^k\to\mathbb
R^n$ be a $C^2$-function defined on an open set $U$ and
$F:\Omega\subset\mathbb R^n\to\mathbb R^n$ a $\mathbb
A$-differentiable vector field defined on an open set $\Omega$
with $\varphi(U)\subset\Omega$. For every initial condition
$w_0\in\Omega$ there exists an unique $\varphi\mathbb
A$-differentiable function $w:V_{\tau_0}\subset\mathbb
R^k\to\mathbb R^n$ with $w(\tau_0)=w_0$ and satisfying
(\ref{apode}), where $V_{\tau_0}\subset U$ is a neighborhood of
$\tau_0$. \ethm \noindent\textbf{Proof.} Define
$$
w_{n+1}(\tau)=\int_{\tau_0}^\tau F\circ w_n(v)d\varphi(v'),\qquad
w_0(v)=w_0.
$$
The function $w_0(v)$ is $\varphi\mathbb A$-differentiable, and by
(\ref{ctfc}) we have that $w_1(v)$ is $\varphi\mathbb
A$-differentiable for all $n\in\mathbb N$. Thus, we apply
induction and show that $w_n(v)$ is $\varphi\mathbb
A$-differentiable. The remaining arguments are similar to the
usual Existence and Uniqueness Theorem for ordinary differential
equations. $\Box$

Let $\mathbb{A}$ be an algebra which as linear space is
$\mathbb{R}^{n}$, and
$\varphi:V\subset\mathbb{R}^{n}\to\mathbb{R}^{n}$ a differentiable
function defined on open set $V$. Consider a function
$F:\Omega\subset\mathbb{R}^{k}\times\mathbb{R}^{n}\to\mathbb{A}$
defined on an open set $\Omega$. We say $F$ is
\emph{$(\varphi\mathbb{A},\mathbb{A})$-differentiable} if
$F(\tau,w)$ as a function of $\tau$ (with $w$ being fixed) is
$\varphi\mathbb{A}$-differentiable and as a function of $w$ (with
$\tau$ being fixed) is $\mathbb{A}$-differentiable. For the
identity map $\varphi:\mathbb{A}\to\mathbb{A}$, we also say $F$ is
\emph{$(\mathbb{A},\mathbb{A})$-differentiable} if $F(\tau,w)$ is
$(\varphi\mathbb{A},\mathbb{A})$-differentiable. We say a function
$f:U\subset\mathbb R^{n+1}\to\mathbb R^{n}$ defined on an open set
$U$ has a \emph{lifting}
$F:\Omega\subset\mathbb{A}\times\mathbb{A}\to\mathbb{A}$ if
$f(t,x)=F(te,x)$ for all $(t,x)\in U$.

A non-autonomous $\varphi\mathbb A$-ordinary differential equation
($\varphi\mathbb A$-ODE) is written by
\begin{equation}\label{varphi:ode}
    w'_\varphi=F(\tau,w),\qquad w(\tau_0)=w_0,
\end{equation}
where finding a solution is understood as the problem of finding a
$\varphi\mathbb A$-differentiable function
$w:V_{\tau_0}\subset\mathbb R^k\to\mathbb R^n$ defined in a
neighborhood $V_{\tau_0}$ of $\tau_0$ whose $\varphi\mathbb
A$-derivative $w'_\varphi(\tau)$ satisfies
\begin{equation}\label{varph:sol}
dw_\tau=F(\tau,w(\tau))d\varphi_\tau.
\end{equation}
The corresponding existence and uniqueness of solutions can be
stated for $(\varphi\mathbb{A},\mathbb{A})$-differentiable
functions $F=F(\tau,w)$.

Given an algebra $\mathbb A$, we say a function
$H:\Omega\subset\mathbb A\times\mathbb A\to\mathbb A$ is of
\emph{$\mathbb A$-separable variables} if $H(\tau,w)=K(\tau)L(w)$
for certain functions $K$ and $L$ which we call \emph{$\mathbb
A$-factors} of $H$.

The $\mathbb{A}$-line integral is defined by the
$\varphi\mathbb{A}$-line integral when
$\varphi:\mathbb{R}^{n}\to\mathbb{R}^{n}$ is the identity map.
Some $\varphi\mathbb{A}$-differential equations can be solved, as
we can see in the following proposition.

\bprop\label{profpin} Suppose $H(\tau,w)=K(\tau)L(w)$ is a
$(\varphi\mathbb{A},\mathbb{A})$-differentiable function of
$\mathbb A$-separable variables, where $L$ has image contained in
the regular set of $\mathbb{A}$. If $w$ is implicitly defined by
\begin{equation}\label{vplip}
    \int^{w}\frac{dv}{L(v)}=\int^{\tau}K(s)d\varphi(s'),
\end{equation}
where the left hand of (\ref{vplip}) denotes the $\mathbb A$-line
integral and the right hand of (\ref{vplip}) denotes the
$\varphi\mathbb A$-line integral, then $w$ is a $\varphi\mathbb
A$-differentiable function of $\tau$ which solves the
$\varphi\mathbb{A}$-differential equation (\ref{varphi:ode}).
\eprop \noindent \noindent\textbf{Proof.} If $w$ is implicitly
defined by (\ref{vplip}) as a function of $\tau$, then by applying
Lemma \ref{l1} to the left hand of (\ref{vplip}) we calculate
$$
\left(\int^{w}\frac{dv}{L(v)}\right)_{\varphi}' =
\left(\frac{d}{dw}\int^{w} \frac{dv}{L(v)}\right)w_{\varphi}' =
\frac{w_{\varphi}'}{L(w)}.
$$
Since $K(\tau)$ is $\varphi\mathbb{A}$-differentiable, by
Corollary \ref{ctfc} we have
$$
\left(\int^{\tau}K(s)d\varphi(s')\right)_{\varphi}'=K(\tau).
$$
Therefore, $w(\tau)$ is a solution of the
$\varphi\mathbb{A}$-differential equation (\ref{varphi:ode}).
$\Box$

\bcor Consider an algebra $\mathbb A$, a differentiable function
in the usual sense $\varphi$, and the $\varphi\mathbb{A}$-ODE
$w_{\varphi}' =K(\tau)w^{2}$ where $K$ is a
$\varphi\mathbb{A}$-differentiable function with
$\varphi\mathbb{A}$-antiderivative $H$. Then, by (\ref{vplip})
$-w^{-1}=H(\tau)+C$, where $C$ is a constant in $\mathbb{A}$.
Thus, the solutions are given by
\begin{equation}\label{evode}
    w(\tau)=\frac{-e}{H(\tau)+C}.
\end{equation}
\ecor

\bejem\label{sistcua} Consider $\mathbb A=\mathbb
A^{2}_1(\alpha,\beta)$ (see Section \ref{apvf}) and a
differentiable function $\varphi$. The $\varphi\mathbb{A}$-ODE
\begin{eqnarray*}
  (xe_1+ye_2)'_\varphi &=&(k(t,s)(x^{2}+\alpha y^{2})+\alpha l(t,s)(2xy+\beta
y^{2}))e_1 \\
   &+& (k(t,s)(2xy+\beta y^{2})+l(t,s)(x^{2}+\alpha y^{2})+\beta
l(t,s)(2xy+\beta y^{2}))e_2
\end{eqnarray*}
can be written by
$$
w'_\varphi=K(\tau)w^{2},
$$
where $K=(k,l)$ is a $\varphi\mathbb{A}$-differentiable function
of $\tau=(t,s)$. \eejem

\bcor Suppose $H(\tau,w)=\varphi(\tau)L(w)$, where $\varphi$ is
differentiable in the usual sense and $L$ is an
$\mathbb{A}$-differentiable function with image contained in the
regular set of $\mathbb{A}$, thus $H$ is
$(\varphi\mathbb{A},\mathbb{A})$-differentiable. Furthermore,
$\varphi\mathbb A$-line integral of $K(\tau)$ satisfies
$$
\int^{\tau}\varphi(s)d\varphi(s')=\frac{(\varphi(\tau))^{2}}{2}.
$$
If $w$ is implicitly defined by
$$
    \int^{w}\frac{dv}{L(v)}=\frac{(\varphi(\tau))^{2}}{2},
$$
then $w$ is a $\varphi\mathbb A$-differentiable function of $\tau$
which solves the $\varphi\mathbb{A}$-differential equation
(\ref{varphi:ode}), where the left hand of (\ref{vplip}) denotes
the $\mathbb A$-line integral. \ecor

We consider the following example.

\bejem\label{exam:varphiodef} Consider the $\varphi\mathbb{A}$-ODE
$w_{\varphi}' =K(\tau)$, where $K$ is a differentiable function in
the usual sense, $\varphi=K$, and $\mathbb A$ is an algebra. Then,
by (\ref{vplip}) $w=\frac{(K(\tau))^{2}}{2}+C$, where $C$ is a
constant in $\mathbb{A}$. Thus, the solutions are given by
$$
w(\tau)=\frac{(K(\tau))^{2}}{2}+C.
$$
\eejem

Another example can be given by

\bejem\label{exam:varphiode} Consider an algebra $\mathbb A$, a
differentiable function in the usual sense $K$, and the
$\varphi\mathbb{A}$-ODE $w_{\varphi}' =K(\tau)w^{2}$ for
$\varphi=K$. Then, by (\ref{vplip})
$-w^{-1}=\frac{(K(\tau))^{2}}{2}+C$, where $C$ is a constant in
$\mathbb{A}$. Thus, the solutions are given by
$$
w(\tau)=\frac{-e}{\frac{(K(\tau))^{2}}{2}+C}.
$$
\eejem

If $k=n$ and $F$ is a $\varphi\mathbb A$-differentiable vector
field, we also say $F$ is a $\varphi\mathbb A$-\emph{algebrizable
vector field}. Consider an ordinary differential equation (ODE)
\begin{equation}\label{ode}
    \frac{dx}{dt}=f(t,x),
\end{equation}
where $f:U\subset\mathbb R^{n+1}\to\mathbb R^{n}$ and $U$ is an
open set. Suppose $\mathbb{A}$ is an algebra such that $f$ has a
lifting $F:\Omega\subset\mathbb{A}\times\mathbb{A}\to\mathbb{A}$
for some set $\Omega$, then an $\mathbb{A}$-ODE associated with
(\ref{ode}) is
\begin{equation}\label{aode}
    \frac{dw}{d\tau}=F(\tau,w),
\end{equation}
where $\frac{dw}{d\tau}$ denotes the $\mathbb{A}$-derivative. When
$F$ is $(\mathbb{A},\mathbb{A})$-differentiable, solutions
$w(\tau)$ of (\ref{aode}) define solutions $x(t)=w(te)$ of
(\ref{ode}), where $e$ is the unit of $\mathbb{A}$. If $H$ is the
function defined by $H(\tau,w)=(d\varphi_\tau(e))^{-1} F(\tau,w)$,
the $\varphi\mathbb A$-differential equation associated with
(\ref{ode}) is given by
\begin{equation}\label{ode:varphi}
    w_\varphi'=H(\tau,w).
\end{equation}
When $H$ is $(\varphi\mathbb{A},\mathbb{A})$-differentiable,
solutions $w(\tau)$ of (\ref{ode:varphi}) define solutions
$x(t)=w(te)$ of (\ref{ode}).

We say: a) the ODE (\ref{ode}) is $\mathbb{A}$-\emph{algebrizable}
if $F$ is $(\mathbb{A},\mathbb{A})$-differentiable, and b) the ODE
(\ref{ode}) is $\varphi\mathbb{A}$-\emph{algebrizable} if $H$ is
$(\varphi\mathbb{A},\mathbb{A})$-differentiable.

The $\varphi\mathbb{A}$-algebrizability for some differentiable
function $\varphi$ of a given ODE as (\ref{ode}) can be
investigated in a similar way as it is done in \cite{AMEENADE}.

\bprop Consider an algebra $\mathbb A$, a differentiable function
$\varphi$ in the usual sense, a $\varphi\mathbb{A}$-differentiable
function $K$, and $g(t)=dH_{(t,0)}(e_1)$, where $H$ is a
$\varphi\mathbb{A}$-antiderivative of $K$. Under these conditions
the $\varphi\mathbb A$-ODE $w'_\varphi=K(\tau)w^{2}$ is a
$\varphi\mathbb{A}$-\emph{algebrization} of the ODE
\begin{equation}\label{spa}
    \frac{dw}{dt}=g(t)w^{2}.
\end{equation}
Thus, by Proposition \ref{profpin} solutions of (\ref{spa}) are
given by
\begin{equation}\label{evode}
    (x(t),y(t))=\frac{-e}{H(t,0)+C}.
\end{equation}
\eprop \noindent\textbf{Proof.} Consider a solution $w$ of the
$\varphi\mathbb A$-ODE $w'_\varphi=K(\tau)w^{2}$, then
\begin{eqnarray*}
  \frac{dw(t,0)}{dt} &=& dw_{(t,0)}(e_1)=w'_\varphi(t,0) d\varphi_{(t,0)}(e_1) \\
    &=& K(t,0)w^{2}d\varphi_{(t,0)}(e_1)\\
    &=& dH_{(t,0)}(e_1)w^{2}.\,\,\,\,\Box
\end{eqnarray*}

\subsection{A $\varphi\mathbb A$-algebrizable vector field
associated to triangular billiards}\label{vaavfatb}

Consider the four-dimensional vector field defined on $\mathbb
R^4$
\begin{equation}\label{cvrd4}
    f(x_1,y_1,x_2,y_2)=\left(%
\begin{array}{c}
  b(x_1^2-y_1^2)-(b+c)(x_1x_2-y_1y_2) \\
  2bx_1y_1-(b+c)(x_1y_2+x_2y_1) \\
  a(x_2^2-y_2^2)-(a+c)(x_1x_2-y_1y_2) \\
  2ax_2y_2-(a+c)(x_1y_2+x_2y_1) \\
\end{array}%
\right)^T.
\end{equation}
It can be written as the two dimensional complex vector field
\begin{equation}\label{cvf1}
    F(u,v)=(bu^2-(b+c)uv,av^2-(a+c)uv),
\end{equation}
for $u=(x_1,y_1)$ and $v=(x_2,y_2)$.

If $\mathbb A$ is the linear space $\mathbb C^2$ endowed with an
algebra structure over $\mathbb C$. We still using $\{e_1,e_2\}$
for the standard basis of $\mathbb C^2$ as a linear space over
$\mathbb C$. For a $\mathbb C$-linear transformation
$\varphi:\mathbb C^2\to\mathbb C^2$, the $\varphi\mathbb
A$-differentiability, the $\varphi\mathbb A$-differential
equations, and the $\varphi\mathbb A$-algebrizability of vector
fields and autonomous ordinary differential equations are defined
in the same way for algebras over $\mathbb C$ as the definitions
for algebras over the real field $\mathbb R$, given in Section
\ref{spadr} and Subsection \ref{spade}. Let
$G:\Omega\subset\mathbb C^2\to\mathbb C^2$ be a complex vector
field defined on an open set $\Omega$. A $\varphi\mathbb
A$-\emph{differential equation} is
\begin{equation}\label{apode}
    w'_\varphi=G(\tau,w),\qquad w(\tau_0)=w_0,
\end{equation}
which is understand as the problem of finding a $\varphi\mathbb
A$-differentiable function $w:V_{\tau_0}\subset\mathbb
C^2\to\mathbb C^2$ defined in a neighborhood
$V_{\tau_0}\subset\Omega$ of $\tau_0$ such that
$w'_\varphi(\tau)=G(\tau,w(\tau))$ for all $\tau\in V_{\tau_0}$.

By applying Corollary \ref{ciaab} to $F$ given in (\ref{cvf1})
$$
M(b_i,\alpha,\beta)=\left(%
\begin{array}{cccc}
  2b & 0 & -\beta(a+c)-(b+c) & a+c \\
  -\beta(a+c)-(b+c) & a+c & 2\beta a & -2a \\
  0 & -2b & -\alpha(a+c) & b+c \\
  -\alpha(a+c) & b+c & 2\alpha a & 0 \\
\end{array}%
\right).
$$
For
\begin{equation}\label{ab}
    \alpha=-\frac{(b+c)^2}{(a+c)^2},\qquad\beta=-2\frac{b+c}{a+c}+\frac{4ab}{(a+c)^2},
\end{equation}
we have $\det(M(b_i,\alpha,\beta))=0$ and
$$
v=\left(1,-\frac{b+c}{a+c},0,-\frac{2b}{a+c}\right)
$$
is in the orthogonal complement of the files of
$M(b_i,\alpha,\beta)$. Thus, by Corollary \ref{ciaab} $F$ is
$\varphi\mathbb A$-differentiable for $\mathbb A=\mathbb
A^{2}_1(\alpha,\beta)$ given in (\ref{pa1ab}) and
\begin{equation}\label{eq:vp}
    \varphi(u,v)=\left(u-\frac{(b+c)}{2b}v,-\frac{(a+c)}{2b}v\right).
\end{equation}

Now, we will look for the expression of $F$ in terms of the
variable $w$ of $\mathbb A$. By evaluating $F$ along the unit
$e\in\mathbb A$ we have
$$
 F(u,0) = b(u,0)^2 = b\left(\varphi(u,0)\right)^2 \\
   = b(\varphi(u,0))^2.
$$
Thus, the complex vector field $F$ in terms of the variable
$w=(u,v)$ of $\mathbb A$ is given by
\begin{equation}\label{edcva}
    F(w)=b(\varphi(w))^2,
\end{equation}
where $(\varphi(w))^2$ is defined by the $\mathbb A$-product.
Therefore, $F$ is $\varphi\mathbb A$-algebrizable. The
corresponding complex ordinary differential equation can be
written by
\begin{equation}\label{ccode}
    \frac{dw}{dz}=b(\varphi(w))^2,
\end{equation}
and its associated $\mathbb{A}$-ODE is given by
\begin{equation}\label{acodebillar}
\frac{dw}{d\tau}=b(\varphi(w))^2.
\end{equation}
Since $d\varphi_w=\varphi$ and $\varphi(e)=e$, the corresponding
$\varphi\mathbb{A}$-ODE is given by
\begin{equation}\label{vphidebillar}
w'_\varphi=b(\varphi(w))^2.
\end{equation}

\subsection{Canonical coordinates of vector fields}

Consider the vector field $F$ given in Section (\ref{cvf1}). A
diffeomorphism $R=(s,r)$ is said to define \emph{canonical
coordinates} for $F$ if $dR(F)=(1,0)$.

Let $F$ be a vector field, $\mathbb A$ an algebra, and $R$ a
$\varphi\mathbb A$-differentiable function which defines canonical
coordinates for $F$, for some differentiable function $\varphi$.
Since $dR(F)=R'_\varphi d\varphi(F)=e$, where $e$ is the unit of
$\mathbb A$, we have $R'_\varphi =(d\varphi(F))^{-1}$. Thus, the
$\varphi\mathbb A$-ODE
\begin{eqnarray*}
  R'_\varphi &=& \left[\left(%
\begin{array}{cc}
  \varphi_{1x} & \varphi_{1y} \\
  \varphi_{2x} & \varphi_{2y} \\
\end{array}%
\right)\left(%
\begin{array}{c}
  f \\
  g \\
\end{array}
\right)\right]^{-1} \\
    &=& \left(%
\begin{array}{c}
 \varphi_{1x}f+\varphi_{1y}g \\
 \varphi_{2x}f+\varphi_{2y}g \\
\end{array}%
\right)^{-1}.
\end{eqnarray*}

If $d\varphi(F)$ is $\varphi\mathbb A$-differentiable, invertible
with respect to $\mathbb A$, and $H$ is the $\varphi\mathbb
A$-antiderivative of $(d\varphi(F))^{-1}$, then the solutions of
$R'_\varphi=(d\varphi(F))^{-1}$ are given by $R(z)=H(z)+C$, where
$C$ is a constant in $\mathbb A$. In this case, by Lemma \ref{l2}
and properties of $\varphi\mathbb A$-differentiability we have
$F=d\varphi^{-1}\circ G\circ\varphi$, where $G$ is $\mathbb
A$-differentiable.

Therefore, the vector fields which can be solved in these way are
basically the $\mathbb A$-differentiable vector fields.

\subsection{On $\varphi\mathbb M$-differentiability for matrix algebras $\mathbb M$}

Let $\mathbb M$ be a commutative matrix algebra with base
$\beta=\{R_1,R_2,\cdots,R_l\}$. Consider two differentiable
functions $f,\varphi:U\subset\mathbb R^{k}\to\mathbb M$ in the
usual sense, where $U$ is an open set. We say $f$ is a
$\varphi\mathbb M$-\emph{differentiable function} if there exists
a function $f'_\varphi:U\subset\mathbb R^{k}\to\mathbb M$ such
that the usual differential $df$ satisfies $df_u=f'_\varphi(u)
d\varphi_u$. That is, $df_u(v)=f'_\varphi(u) d\varphi_u(v)$  for
all $v\in\mathbb R^{k}$, where $f'_\varphi(u) d\varphi_u(v)$
denotes the product of the matrices $f'_\varphi(u)$ and
$d\varphi_u(v)$.

For this definition, polynomial functions like
(\ref{Apolynomial:function}), rational functions in the same way,
and exponential function of the form $u\mapsto e^{f(u)}$, where
$f$ is $\varphi\mathbb M$-differentiable, are well defined and
they satisfy the usual rules of differentiation, as we comment in
Section \ref{varphi:derivative}. Thus, for this kind of derivative
also we can consider $\varphi\mathbb A$-differential equations
\begin{equation}\label{separable:vamodes}
   w'_\varphi=H(\tau,w)
\end{equation}
which are matrix differential equations. If $H(\tau,w)$ is
$(\varphi\mathbb M,\mathbb M)$-differentiable (see Section
\ref{spade}) and $H(\tau,w)=K(\tau)L(w)$, then we have Proposition
\ref{vplip} and the results given for $\varphi\mathbb
A$-derivatives and $\varphi\mathbb A$-differential equations. For
example,
\begin{equation}\label{exp:dif:eq}
    w'_\varphi=w
\end{equation}
has the unique solution $w(\tau)=Ce^{\varphi(\tau)}$ with
$w(\tau_0)=Ce^{\varphi(\tau_0)}$, where $C\in\mathbb M$.

This definition of $\varphi\mathbb M$-differentiability allows us
to have matrix differential equations, and it can be used for
giving solutions of some systems of partial differential
equations, it is done in \cite{CAL2}. In the following section we
show solutions of systems of PDEs which can be obtained in this
way.

\section{On solutions of PDEs}\label{sol:singles:pdes}

\subsection{A class of first order PEDs}

We consider a PDE of the form
\begin{equation}\label{firstorder:linear:pde}
    au_x+bv_x-cu_y-dv_y=0,
\end{equation}
where $a$, $b$, $c$, and $d$ are functions of $x$ and $y$. This is
an homogeneous first order partial differential equation with two
dependent variables $u$, $v$ and two independent variables $x$,
$y$.

If $a$, $b$, $c$ and $d$ are constants, equation
(\ref{firstorder:linear:pde}) can be written by
\begin{equation}\label{pde:accommodated}
    (au+bv)_x=(cu+dv)_y.
\end{equation}
A standard technique in the research for solutions of differential
equations is the separation of variables, that in this case it can
be applied in the way
\begin{equation}\label{usual:asumption}
    (au+bv)_x=f'(x)g'(y)=(cu+dv)_y,
\end{equation}
then
$$
    au+bv = f(x)g'(y)+h_1(y), \quad cu+dv =
    f'(x)g(y)+h_2(x).
$$
Thus, if $ad-bc\neq 0$,
\begin{eqnarray*}
  u(x,y) &=& \frac{1}{ad-bc}(df(x)g'(y)-bf'(x)g(y)+dh_1(y)-bh_2(x)),\\
  v(x,y) &=& \frac{1}{ad-bc}(af'(x)g(y)-cf(x)g'(y)+ah_2(x)-ch_1(y)).
\end{eqnarray*}

The aim of this section is to avoid assumption
(\ref{usual:asumption}) and to propose families of functions which
satisfy (\ref{pde:accommodated}). We obtain families of solutions
of (\ref{firstorder:linear:pde}) which are families of
$\varphi\mathbb A$-differentiable functions for adequate functions
$\varphi$ and algebras $\mathbb A$.

\subsection{On $\varphi\mathbb A$-CREs for $\mathbb A=\mathbb A^{2}_1(\alpha,\beta)$,
$\mathbb A=\mathbb A^{2}_2(\alpha,\beta)$, and $\mathbb A=\mathbb A^{2}_{1,1}$}

Consider the $\varphi\mathbb A$-CREs given in (\ref{gcre}). For
$\mathbb A=\mathbb A^{2}_1(\alpha,\beta)$ defined by
(\ref{ffra1ab}), they are
\begin{equation}\label{varphi:algebra1:GCREsc}
    \begin{array}{c}
      \varphi_{1y}r_x+\alpha\varphi_{2y}s_x-\varphi_{1x}r_y-\alpha\varphi_{2x}s_y=0, \\
      \varphi_{2y}r_x+(\varphi_{1y}+\beta\varphi_{2y})s_x-\varphi_{2x}r_y-(\varphi_{1x}+\beta\varphi_{2x})s_y=0. \\
    \end{array}
\end{equation}
For $\mathbb A=\mathbb A^{2}_2(\alpha,\beta)$ defined by
(\ref{ffra2ab}), they are
\begin{equation}\label{varphi:algebra2:GCREsc}
    \begin{array}{c}
      (\alpha\varphi_{1y}+\varphi_{2y})r_x+\varphi_{1y}s_x-(\alpha\varphi_{1x}+\varphi_{2x})r_y-\varphi_{1x}s_y=0, \\
      \beta\varphi_{1y}r_x+\varphi_{2y}s_x-\beta\varphi_{1x}r_y-\varphi_{2x}s_y=0. \\
    \end{array}
\end{equation}
For  $\mathbb A=\mathbb A^{2}_{1,1}$ defined by (\ref{ffra3}),
they are
\begin{equation}\label{varphi:algebra3:GCREsc}
     \varphi_{1y}r_x-\varphi_{1x}r_y=0,\qquad \varphi_{2y}s_x-\varphi_{2x}s_y=0.
\end{equation}

\subsection{Solutions defined by $\varphi\mathbb A$-differentiable
functions}

In the following Theorem we give conditions under which for some
algebra $\mathbb A$ of the type $\mathbb A^{2}_1(\alpha,\beta)$,
$\mathbb A^{2}_2(\alpha,\beta)$, or $\mathbb A^{2}_{1,1}$, which
are defined in Section \ref{apvf}, and for a differentiable
function $\varphi$, the $\varphi\mathbb A$-differentiable
functions are solutions of (\ref{firstorder:linear:pde}).

\bthm Let $\varphi=(\varphi_1,\varphi_2)$ be a differentiable
function in the usual sense. The $\varphi\mathbb A$-differentiable
functions are solutions of the system of homogeneous first order
linear partial differential equation with two dependent variables
$r$ and $s$ given by
\begin{small}
\begin{equation}\label{lineareq:canbesolved:A11}
    (k\varphi_{1y}+l\varphi_{2y})r_x+(l\varphi_{1y}+(k\alpha+l\beta)\varphi_{2y})s_x
    -(k\varphi_{1x}+l\varphi_{2x})r_y-(l\varphi_{1x}+(k\alpha+l\beta)\varphi_{2x})s_y=0,
\end{equation}
\begin{equation}\label{lineareq:canbesolved:A12}
    (k\varphi_{2y}+(k\alpha+l\beta)\varphi_{1y})r_x+(k\varphi_{1y}+l\varphi_{2y})s_x
    -((k\alpha+l\beta)\varphi_{1x}+k\varphi_{2x})r_y-(k\varphi_{1x}+l\varphi_{2x})s_y=0,
\end{equation}
\begin{equation}\label{lineareq:canbesolved:A13}
    k\varphi_{1y}r_x+l\varphi_{2y}s_x-k\varphi_{1x}r_y-l\varphi_{2x}s_y=0,
\end{equation}
\end{small}
for $\mathbb A=\mathbb A^{2}_1(\alpha,\beta)$, $\mathbb A=\mathbb
A^{2}_2(\alpha,\beta)$, and $\mathbb A=\mathbb A^{2}_{1,1}$,
respectively, where $k$ and $l$ are functions of $(x,y)$.
\ethm\bdemo Equation (\ref{lineareq:canbesolved:A11}) is equal to
$k$ by first equation of (\ref{varphi:algebra1:GCREsc}) plus $l$
by the second. Equations (\ref{lineareq:canbesolved:A12}) and
(\ref{lineareq:canbesolved:A13}) can be obtained in a similar way
from (\ref{varphi:algebra2:GCREsc}) and
(\ref{varphi:algebra1:GCREsc}), respectively. \edemo

\bcor If there exists a differentiable function
$\varphi=(\varphi_1,\varphi_2)$ in the usual sense and constants
$\alpha$, $\beta$, $k$, and $l$ such that
\begin{eqnarray*}
  a &=& k\varphi_{1y}+l\varphi_{2y}, \\
  b &=& l\varphi_{1y}+(k\alpha+l\beta)\varphi_{2y}, \\
  c &=& k\varphi_{1x}+l\varphi_{2x}, \\
  d &=& l\varphi_{1x}+(k\alpha+l\beta)\varphi_{2x},
\end{eqnarray*}
then the $\varphi\mathbb A$-differentiable functions are solutions
of (\ref{firstorder:linear:pde}) for $\mathbb A=\mathbb
A^{2}_1(\alpha,\beta)$. \ecor

\bcor\label{example:first} Suppose the $\alpha+\beta\neq 1$,
$\mathbb A$ is the algebra $\mathbb A^{2}_1(\alpha,\beta)$, and
$$
\varphi(x,y)=\left(\frac{-c(\alpha+\beta)+d}{\alpha+\beta-1}x+\frac{a(\alpha+\beta)-b}{\alpha+\beta-1}y,\frac{-c+d}{\alpha+\beta-1}x+\frac{a-b}{\alpha+\beta-1}y\right),
$$
then the $\varphi\mathbb A$-differentiable functions
$w(x,y)=(u(x,y),v(x,y))$ are solutions of
(\ref{firstorder:linear:pde}).\ecor

\bejem Suppose situation of Corollary \ref{example:first} with
$\alpha=0$ and $\beta=0$. That is, $\mathbb A$ is the algebra
$\mathbb A^{2}_1(0,0)$ and the $\mathbb A$-product is determined
by $e_1e_1=e_1$, $e_1e_2=e_2$, $e_2e_2=0$. Thus, we have
$$
\varphi(x,y)=(dx+by,(c-d)x+(a-b)y).
$$
Since $\varphi$ is a $\varphi\mathbb A$-differentiable function we
have that the pair of functions
\begin{eqnarray*}
  u(x,y) &=& dx+by, \\
  v(x,y) &=& (c-d)x+(a-b)y,
\end{eqnarray*}
is a solution of (\ref{firstorder:linear:pde}). The function
$\varphi^{2}$ is $\varphi\mathbb A$-differentiable, in this case
the solution $(u,v)$ of (\ref{firstorder:linear:pde}) is given by
\begin{eqnarray*}
  u(x,y) &=& d^{2}x^{2}+2bdxy+b^{2}y^{2}, \\
  v(x,y) &=& 2(cd-d^{2})x^{2}+2(ad+bc-2bd)xy+2(ab-b^{2})y^{2}.
\end{eqnarray*}

All polynomial, rational, and exponential functions of
$\varphi(x,y)$ with respect to $\mathbb A$ are $\varphi\mathbb
A$-differentiable. Therefore, all of these functions are solutions
of (\ref{firstorder:linear:pde}). \eejem

By Corollary \ref{example:first} we found families of solutions of
(\ref{firstorder:linear:pde}) for constant coefficients $a$, $b$,
$c$, and $d$. If we set $\alpha=-1$ and $\beta=0$, then the
solutions can be of the form $f\circ \varphi$ for all $f$
differentiable functions in the complex sense.

\subsection{EDPs with more dependent and independent
variables}\label{s4}

If we consider the PDE
\begin{equation}\label{firstorder:linear:pde3}
    a_1u_x+b_2v_x+c_1w_x+a_2u_y+b_2v_y+c_2w_y+a_3u_z+b_2v_z+c_2w_z=0,
\end{equation}
we may use the three-dimensional algebras given in \cite{Dyg3}, to
obtain solutions of this type of PDEs. Therefore, if
(\ref{firstorder:linear:pde3}) has constant coefficient, then
there exist three dimensional algebras $\mathbb A$ such that the
$\varphi\mathbb A$-differentiable functions are solutions of
(\ref{firstorder:linear:pde3}).

\subsection{On algebrizability of systems of two
PDEs}\label{algsection}

The algebrizability systems of two first order partial
differential equations with two dependent variables and two
independent variables is considered in \cite{CAL2}, where
solutions of system are constructed
\begin{align}\label{eqn:pde:ex1}
a_1\frac{\partial y}{\partial x}+ \frac{\partial y}{\partial t}=
-b_1y+b_1z, \\ \notag -a_2\frac{\partial z}{\partial x}+
\frac{\partial z}{\partial t} = b_2y-b_2z,
\end{align}
with constant coefficients which appears in \cite{MSEQandA}.
Examples of these solutions are
$$
\begin{array}{c}
  y(x,t)= c_1e^{h_1(x,t)}\cos h_2(x,t)+c_2e^{h_1(x,t)}\sin h_2(x,t), \\
  z(x,t)= -c_1e^{h_1(x,t)}\sin h_2(x,t)+c_2e^{h_1(x,t)}\cos h_2(x,t),\\
\end{array}
$$
for all $c_1,c_2\in\mathbb R$, where
$$
h_1(x,t)=\frac{-b_1+b_2}{a_1+a_2}x+\frac{-a_1b_2-a_2b_1}{a_1+a_2}t,\qquad
h_2(x,t)=\frac{b_1+b_2}{a_1+a_2}x+\frac{-a_1b_2+a_2b_1}{a_1+a_2}t.
$$
Another solutions are
$$
\begin{array}{c}
  y(x,t)= c_1e^{-h(x,t)}\cosh h(x,t)+c_2e^{-h(x,t)}\sinh h(x,t), \\
  z(x,t)= c_1e^{-h(x,t)}\sinh h(x,t)+c_2e^{-h(x,t)}\cosh h(x,t),\\
\end{array}
$$
for all $c_1,c_2\in\mathbb R$, where
$$
h(x,t)=\frac{b_1-b_2}{a_1+a_2}x+\frac{a_1b_2+a_2b_1}{a_1+a_2}t.
$$

\section{Acknowledgements}
The author thanks James S. Cook for the suggestions given that
have helped to improve this work.

\end{document}